%% file: one.tex
\documentclass{article}

\input{diagrams}
\diagramstyle[repositionpullbacks,midshaft,PostScript=dvips,nohug,scriptlabels,size=2em] 
\usepackage{graphics}
\usepackage{color}
\usepackage[mathcal,mathscr]{euscript}

\input{abstd}

\input{abmore}

\input{thmlist}

\newcommand{\xrefstyle}[1]{#1}
\newcommand{\secthecantorset}{\xrefstyle{4}}
\newcommand{\lemmafixedpointcomponents}{\xrefstyle{1.1}}
\newcommand{\egrecFreyd}{\xrefstyle{2.1}}
\newcommand{\egrecSierpinski}{\xrefstyle{2.10}}
\newcommand{\egrecifs}{\xrefstyle{2.11}}
\newcommand{\egrecbary}{\xrefstyle{2.12}}

\begin{document}

\sloppy

\include{front1}

\input{intro1}

\input{ssss}

\input{univsoln}

\appendix

\input{solvability}

\input{refs1}
\end{document}

%% file: abstd.tex
\usepackage{amssymb}
\newcommand{\emptybk}{\:\:}

\newcommand{\dashbk}{-}
\newcommand{\diagspace}{\mbox{\hspace{2em}}}
\newcommand{\bref}[1]{(\ref{#1})}
\newcommand{\ucontents}[2]{\addcontentsline{toc}{#1}{\numberline{}{#2}}}
\newcommand{\mb}[1]{\mathbf{#1}}
\newcommand{\mc}[1]{\mathcal{#1}}
\newcommand{\mr}[1]{\mathrm{#1}}
\newcommand{\cat}[1]{\mc{#1}}
\newcommand{\fcat}[1]{\mb{#1}}
\newcommand{\ovln}[1]{\overline{#1}}
\newcommand{\twid}[1]{\widetilde{#1}}
\newcommand{\slsh}{/\linebreak[0]}
\newcommand{\dt}{.\linebreak[0]}
\newcommand{\such}{\:|\:}

\newcommand{\blob}{{\raisebox{.3ex}{\ensuremath{\scriptscriptstyle{\bullet}}}}}

\newcommand{\implies}{\,\Rightarrow\,}
\newcommand{\op}{\mr{op}}
\newcommand{\id}{\mr{id}}
\newcommand{\Cat}{\fcat{Cat}}
\newcommand{\One}{\fcat{1}}
\newcommand{\Set}{\fcat{Set}}
\newcommand{\Top}{\fcat{Top}}
\newcommand{\ftrcat}[2]{[#1,#2]}
\newcommand{\goesto}{\,\longmapsto\,}
\newcommand{\goiso}{\goby{\diso}}
\newcommand{\og}{\lTo}
\newcommand{\ogby}[1]{\lTo^{#1}}
 
\newcommand{\parpair}[2]{\pile{\rTo^{\scriptstyle #1}\\ 
\rTo_{\scriptstyle #2}}}
\newcommand{\parpairu}{\pile{\rTo\\ \rTo}}
\newarrow{Get}....>
\newarrow{Goesto}|---{->}
\newarrow{Incl}C--->
\newarrow{Mod}--+->
\newcommand{\diso}{\sim}

%% file: abmore.tex
\newcommand{\reals}{\mathbb{R}}
\newcommand{\latin}[1]{\textit{#1}}
\newcommand{\vs}{vs}			
{\begin{array}{l}}%
{\end{array}}
\newcommand{\demph}[1]{\textbf{\textup{#1}}}
\newcommand{\done}{\hfill\ensuremath{\Box}}
\newenvironment{prooflike}[1]{\begin{trivlist}\item\textbf{#1}\ }
{\end{trivlist}}
\newenvironment{proof}{\begin{prooflike}{Proof}}{\end{prooflike}}
\newcommand{\scat}[1]{\mathbb{#1}}
\newcommand{\go}{\rTo\linebreak[0]}
\newcommand{\goby}[1]{\rTo^{#1}\linebreak[0]}
\newcommand{\iso}{\cong}
\newcommand{\nat}{\mathbb{N}}	
\newcommand{\url}[1]{#1}
\newcommand{\pshf}[1]{\ftrcat{#1^\op}{\Set}}
\newcommand{\elt}[1]{\scat{E}\left(#1\right)}
\newcommand{\gomod}{\rMod\linebreak[0]}
\newcommand{\gobymod}[1]{\rMod^{#1}\linebreak[0]}
\newcommand{\ndSet}[1]{[#1, \Set]_{\mr{nondegen}}}
\newcommand{\ndTop}[1]{[#1, \Top]_{\mr{nondegen}}}
\newcommand{\littleproduct}{(\blob\ \ \ \blob)}
\newcommand{\littleequalizer}{(\blob \parpairu \blob)}
\newcommand{\littlepullback}{(\blob \go \blob \og \blob)}
\newcommand{\cstyle}[1]{\textbf{\upshape{#1}}}
\newcommand{\So}{\cstyle{S}}
\newenvironment{condition}{\begin{description}\item[\ \ \ \ \ \,]}{\end{description}}
\newarrow{Qt}---->
\newarrow{Modget}..+.>
\newcommand{\Coalg}[2]{\fcat{Coalg}(#1, #2)}
\newcommand{\frc}[2]{\frac{#1}{#2}}
\newcommand{\half}{\frc{1}{2}}
\newcommand{\dfrc}[2]{\textstyle{\frac{#1}{#2}}}
\newcommand{\dhalf}{\dfrc{1}{2}}
\newcommand{\compt}[1]{{[#1]}}
\newcommand{\bigcompt}[1]{{[\, #1 \,]}}
\newcommand{\ob}{\mr{ob}}		
\newcommand{\pr}{\mr{pr}}
\newcommand{\Res}{\mr{Res}}

\newcommand{\vtr}[1]{#1}
\newcommand{\oba}{\mr{ob}\,}
\newcommand{\blb}{{\scriptscriptstyle{\bullet}}}
\newcommand{\sblb}{\cdot}
\newcommand{\ldlabel}[1]{\makebox[0em]{\hspace*{4.2em}\ensuremath{\scriptstyle#1}}}
\newcommand{\rdlabel}[1]{\makebox[0em]{\hspace*{-4.5em}\ensuremath{\scriptstyle#1}}}

\newcommand{\of}{\,\raisebox{0.08ex}{\ensuremath{\scriptstyle\circ}}\,}
\newcommand{\Cx}{\mathbb{C}}
\newcommand{\sub}{\subseteq}
\newcommand{\catI}{\mathcal{I}}
\newcommand{\ladj}{\dashv}
\newenvironment{diagdiag}
	{\begin{diagram}[size=1.5em]}
	{\end{diagram}}
\newcommand{\cell}[4]{\put(#1,#2){\makebox(0,0)[#3]{\ensuremath{#4}}}}

\definecolor{grey}{gray}{0.5}

\newcommand{\numcell}[3]{\put(#1,#2){\makebox(0,0){\ensuremath{\scriptstyle
#3}}}}

%% file: thmlist.tex
\newtheorem{thm}{Theorem}[section]
\newtheorem{propn}[thm]{Proposition}
\newtheorem{lemma}[thm]{Lemma}
\newtheorem{cor}[thm]{Corollary}
\newtheorem{lotsofremarks}[thm]{Remarks}

\newtheorem{predefn}[thm]{Definition}
\newenvironment{defn}{\begin{predefn}\upshape}{\end{predefn}}
\newtheorem{preexample}[thm]{Example}
\newenvironment{example}{\begin{preexample}\upshape}{\end{preexample}}
\newenvironment{example*}[1]{\begin{preexample}[#1]\upshape}{\end{preexample}}
\newtheorem{prewarning}[thm]{Warning}
\newenvironment{warning}{\begin{prewarning}\upshape}{\end{prewarning}}

%% file: front1.tex
\title{A general theory of self-similarity I}
\author{Tom Leinster}
\date{\normalsize 
Department of Mathematics, University of Glasgow\\
www\dt maths\dt gla\dt ac\dt uk\slsh $\sim$tl\\
tl@$\,\!$maths\dt gla\dt ac\dt uk}

\maketitle
\thispagestyle{empty}

\begin{center}
\textbf{Abstract}
\end{center}
Consider a self-similar space $X$.  A typical situation is that $X$ looks
like several copies of itself glued to several copies of another space $Y$,
and $Y$ looks like several copies of itself glued to several copies of $X$,
or the same kind of thing with more than two spaces.  Thus, the
self-similarity of $X$ is described by a system of simultaneous equations.
Here I formalize this idea and the notion of a `universal solution' of such
a system.  I determine exactly when a system has a universal solution and,
when one does exist, construct it.
\vspace*{2mm}\\
\noindent
A sequel~\cite{SS2} contains further results and examples,
and an introductory article~\cite{GSSO} gives an overview.

\vfill

\tableofcontents

%% file: intro1.tex
\section*{Introduction}
\ucontents{section}{Introduction}

Statements of self-similarity are of two kinds: local and global.  Local
statements say something like `almost any small pattern observed in one
part of the object can be observed throughout the object, at all scales'.
See for instance Chapter~4 of Milnor~\cite{Mil}, where such statements are
made about Julia sets of complex rational functions.  Global statements say
something like `the whole object consists of several smaller copies of
itself glued together'; more generally, there may be a whole family of
objects, each of which can be described as several objects in the family
glued together.  Put another way, the contrast is between bottom-up and
top-down.  This paper and its sequel~\cite{SS2} introduce a theory of
global, or top-down, self-similarity.

Statements of self-similarity can also be divided into the glamorous and
the mundane.  Fractal spaces such as Julia sets, the Cantor set, and
Sierpi\'nski's gasket are certainly self-similar, but there is
self-similarity to be found in more everyday objects.  For instance, any
$n$-manifold is as locally self-similar as could be: every point is locally
isomorphic to every other point.  The closed interval $[0, 1]$ is globally
self-similar, being isomorphic to two copies of itself glued end to end.
We will see that $[0, 1]$ is, in a precise sense, \emph{universal} with
this property (Freyd's Theorem,~\ref{thm:Freyd}).  In the same way,
barycentric subdivision expresses the topological $n$-simplex $\Delta^n$ as
the gluing-together of several smaller copies of itself, and this leads to
a universal characterization of the sequence $(\Delta^n)_{n\geq 0}$ among
all sequences of topological spaces \cite[\egrecbary]{SS2}.

The shape of this paper can be explained by analogy with linear
simultaneous equations.
\begin{enumerate}
\item Consider a sequence $\vtr{x} = (x_1, \ldots, x_n)$ of scalars and a
  system of $n$ simultaneous equations, the $i$th of which expresses $x_i$
  as a linear combination of $x_j$s (that is, $x_i = \sum_j m_{ij} x_j$).
\item The coefficients can be encoded as an $n \times n$ matrix $M$.
\item The system can be written as $\vtr{x} = M\vtr{x}$, and we are
  interested in finding solutions.
\item We consider only nondegenerate solutions $\vtr{x} \neq 0$.
\item We may decide to consider only nondegenerate systems of equations,
  that is, those for which $M$ is nonsingular ($\vtr{x}$ nondegenerate
  implies $M\vtr{x}$ nondegenerate).
\item There are explicit conditions on $M$ equivalent to the existence
  of a nondegenerate solution of the system (for instance, $\det (M - I) =
  0$). 
\item In the case that a nondegenerate solution exists, it can be
  constructed algorithmically.
\end{enumerate}
Here are the analogous steps for self-similarity equations.  
\begin{enumerate}
\item Consider a family $X = (X_a)_{a \in \scat{A}}$ of spaces and a system
  of simultaneous equations, one for each $a \in \scat{A}$, expressing each
  $X_a$ as a gluing-together of $X_b$s~(\S\ref{sec:sss}).
\item The gluing formulas can be encoded as a `two-sided $\scat{A}$-module'
  $M$~(\S\ref{sec:sss}).
\item The system can be written as $X \iso M \otimes X$.  The linear
  analogy is a simplification: there we had scalars, which are either equal
  or not, but here we have spaces, which may be equal, isomorphic, or
  merely connected by a map in one direction or the other.  We are
  interested in `$M$-coalgebras', that is, pairs $(X, \xi)$ where $\xi: X
  \go M \otimes X$.  We are especially interested in $M$-coalgebras $(X,
  \xi)$ possessing a certain universal property.  If this property holds
  then $\xi$ is an isomorphism, so we call such an $(X, \xi)$ a `universal
  solution' of the system (\S\S\ref{sec:sss},~\ref{sec:coalgs}).
\item We consider only `nondegenerate' families $X$
  (\S\S\ref{sec:sss},~\ref{sec:nondegen}). 
\item We consider only nondegenerate systems of equations, that is, those
  for which $M$ is a `nondegenerate' module ($X$ nondegenerate implies $M
  \otimes X$ nondegenerate: \S\S\ref{sec:sss},~\ref{sec:nondegen}).
\item There is an explicit condition on $M$ equivalent to the existence of
  a universal solution of the system (\S\ref{sec:univ-soln} and
  Appendix~\ref{app:solv}; I claim no analogy with determinants).
\item In the case that a universal solution exists, it can be constructed
  algorithmically (\S\S\ref{sec:univ-soln}--\ref{sec:Top-proofs}).
\end{enumerate}
The second paper in this series~\cite{SS2} shows how to recognize universal
solutions, gives examples, classifies those spaces that are self-similar in
at least one way, and uses this classification to reprove some classical
results in topology.

Self-similarity is regarded here as \emph{intrinsic} structure: there is no
ambient space.  (Contrast iterated function systems.)  This is like
considering abstract groups rather than groups of transformations, or
abstract manifolds rather than manifolds embedded in $\reals^n$.

I have been referring to `spaces' with self-similar structure.  Here and
in~\cite{SS2}, the types of space considered are sets and general topological
spaces.  It may be possible to extend the theory to encompass other types of
space, hence other types of self-similarity: conformal, statistical,
type-theoretic (in the sense of computer science), and so on.  The formal
mechanism would be to replace the category of sets, which plays a basic
role in the theory presented here, by a different category of spaces.

Another long-term goal is to develop the algebraic topology of self-similar
spaces, for which the usual homotopical and homological invariants are
often useless: in the case of a connected self-similar subset of the plane, for
example, they only give us $\pi_1$, which is typically either
infinite-dimensional or trivial.  However, a description by a set of
self-similarity equations is discrete and so might provide useful
invariants.

Various other theories are related to this one.  Symbolic
dynamics~\cite{LM} seems most closely related to the case of
\emph{discrete} self-similarity systems, studied in~\S\secthecantorset\
of~\cite{SS2}.  Iterated function systems (\cite{Fal},~\cite{Hut}) are
related, but differ crucially in that they take place inside a fixed
ambient space; see Examples~\egrecSierpinski\ and~\egrecifs\ of~\cite{SS2}
for more.  There is also a paper of Barr~\cite{Barr} with obvious
similarities to the present work.  He discusses terminal coalgebras for an
endofunctor, and the metrics associated with them.  However, the class of
endofunctors he considers is almost disjoint from the class considered
here: the categories on which his endofunctors act always have a terminal
object, and his terminal coalgebras can be constructed as limits
(compare~\ref{warning:fin} below).

\paragraph*{Notation and terminology}  The natural language for this theory
is that of categories.  The following concepts are used; they are all
explained in~\cite{CWM}.

Let $\scat{I}$ be a small category and $\cat{E}$ any category.  A
\demph{diagram of shape $\scat{I}$} (or \demph{over $\scat{I}$}) in
$\cat{E}$ is a functor $D: \scat{I} \go \cat{E}$.  A \demph{cone} on $D$ is
an object $X$ of $\cat{E}$ together with a map $p_i: X \go Di$ for each $i
\in \scat{I}$, such that $(Du) \of p_i = p_j$ for any map $u: i \go j$ in
$\scat{I}$.  When the cone has a certain universal property, it (or
abusively, the object $X$) is called the \demph{limit} of $D$.
\demph{Cocones} and \demph{colimits} are defined dually.

More generally, let $F: \scat{I}^\op \times \scat{I} \go \cat{E}$.  A
\demph{wedge} on $F$ is an object $X$ of $\cat{E}$ together with a map
$p_i: X \go F(i, i)$ for each $i \in \scat{I}$, such that $F(1, u) \of p_i
= F(u, 1) \of p_j$ for any map $u: i \go j$ in $\scat{I}$.  When the wedge
has a certain universal property, it (or $X$) is called the \demph{end} of
$F$, written $X = \int_i F(i, i)$.  \demph{Coends} are defined dually
and written $X = \int^i F(i, i)$.

(Co)limits are a special case of (co)ends: given $D: \scat{I} \go \cat{E}$,
define $F(i, j) = Dj$; then $\int_i F(i, i)$ is the limit of $D$.  So it is
reasonable and convenient to write $\int_i Di$ for the limit of $D$, when
it exists.  Dually, $\int^i Di$ is the colimit of $D$.

The sum (coproduct) of a set-indexed family $(D_i)_{i \in I}$ of objects
is written $\sum_i D_i$.  If $D_i = Y$ for all $i$ then the sum is written
$I \times Y$.  The sum of a finite family $D_1, \ldots, D_n$ of objects is
written $D_1 + \cdots + D_n$, or $0$ if $n = 0$.

Given categories $\cat{A}$ and $\cat{B}$, the category whose objects are
functors from $\cat{A}$ to $\cat{B}$ and whose morphisms are natural
transformations is written $\ftrcat{\cat{A}}{\cat{B}}$.  

A \demph{discrete category} is one in which the only maps are the
identities.  Small discrete categories are therefore just sets.  

$\Top$ is the category of \emph{all} topological spaces and continuous
maps.

The set $\nat$ of natural numbers is taken to include $0$.

\paragraph*{Acknowledgements}
I am very grateful for the encouragement and help of Clemens
Berger, Araceli Bonifant, Edward Crane, Marcelo Fiore, Steve
Lack, John Milnor, Justin Sawon, Carlos Simpson, and Ivan Smith.
This research has relied crucially on the categories mailing list
(see~\cite{Fre}); I thank Bob Rosebrugh, who runs it.  I am
immensely grateful to Jon Nimmo for creating
Figure~\ref{fig:Julia}.  The commutative diagrams were made using
Paul Taylor's macros.

I gratefully acknowledge a Nuffield Foundation award NUF-NAL 04.

%% file: ssss.tex
\section{Self-similarity systems}
\label{sec:sss}

This section concerns two definitions:
\begin{itemize}
\item a `self-similarity system' is a small category $\scat{A}$ together
  with a finite nondegenerate module $M: \scat{A} \gomod \scat{A}$
\item a `universal solution' of a self-similarity system $(\scat{A}, M)$ is
  a terminal $M$-coalgebra.
\end{itemize}
I will explain all this terminology and the underlying ideas using two
examples: a Julia set and a theorem of Freyd.

\begin{figure}
\setlength{\unitlength}{1mm}
\begin{picture}(120,75)
\cell{0}{4}{bl}{\includegraphics{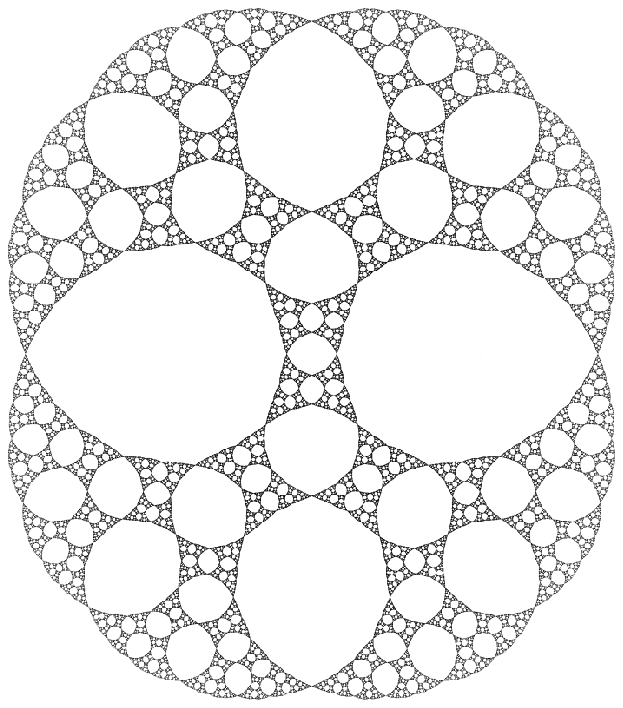}}
\cell{95}{75.5}{t}{\includegraphics{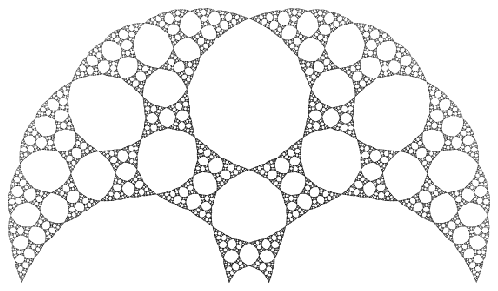}}
\cell{95}{5}{b}{\includegraphics{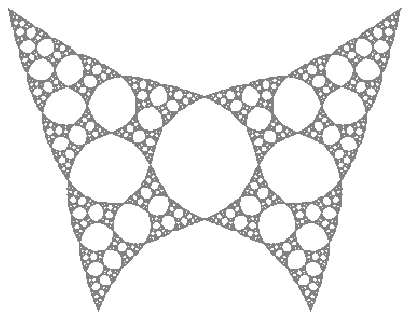}}
\cell{31.5}{0}{b}{\textrm{(a)}}
\cell{95}{39}{b}{\textrm{(b)}}
\cell{72}{46.5}{t}{1}
\cell{92.5}{46.5}{t}{2}
\cell{97.5}{46.5}{t}{3}
\cell{118}{46.5}{t}{4}
\cell{95}{0}{b}{\textrm{(c)}}
\cell{82.5}{4.5}{c}{1}
\cell{73.5}{35}{c}{2}
\cell{116.5}{35}{c}{3}
\cell{107}{4.5}{c}{4}
\end{picture}
\caption{(a) The Julia set of $(2z/(1+z^2))^2$; (b),~(c) two subsets,
  rescaled.  (Image by Jon Nimmo; see also~\cite[Fig.~2]{Mil}
  and~\cite[Fig.~53]{PR})}
\label{fig:Julia}
\end{figure}

Figure~\ref{fig:Julia}(a)
shows a certain closed subset of the Riemann sphere $\Cx \cup \{ \infty
\}$, the Julia set of the function $z \goesto (2z / (1 + z^2))^2$.  (What
follows is informal, and the reader will not need the definition of Julia
set.)  Write $I_1$ for this set, regarded as an abstract topological space.
Evidently $I_1$ has reflectional symmetry in a horizontal axis, so may be
written
\begin{equation}
\label{eq:Julia1}
I_1 = 
\begin{array}{c}
\setlength{\unitlength}{1mm} % eps file is 17 x 20 mm
\begin{picture}(18,20)(-9,-10)
\cell{0.1}{0.1}{c}{\includegraphics{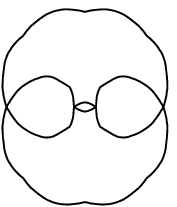}}
\numcell{-7.7}{2}{1}
\numcell{-0.9}{1.8}{2}
\numcell{0.9}{1.8}{3}
\numcell{7.6}{2}{4}
\numcell{-7.7}{-2}{1}
\numcell{-0.9}{-1.8}{2}
\numcell{0.9}{-1.8}{3}
\numcell{7.6}{-2}{4}
\cell{0}{6}{c}{I_2}
\cell{0}{-6}{c}{I_2}
\end{picture}
\end{array}
\end{equation}
where $I_2$ is a certain space with 4 distinguished points, shown in
Figure~\ref{fig:Julia}(b).  In turn, $I_2$ may be regarded as a
gluing-together of subspaces:
\begin{equation}
\label{eq:Julia2}
I_2 = 
\begin{array}{c}
\setlength{\unitlength}{1mm} % eps file is 51 x 30 mm
\begin{picture}(52,30)(-26,0)
\cell{0.2}{0}{b}{\includegraphics{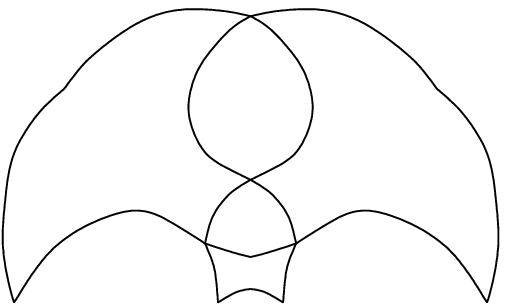}}
\numcell{-23.6}{2.5}{1}
\numcell{-5}{7.8}{2}
\numcell{-2.3}{12.5}{3}
\numcell{-3}{28.5}{4}
\numcell{23.8}{2.5}{1}
\numcell{5.2}{7.8}{2}
\numcell{2.5}{12.5}{3}
\numcell{3}{28.5}{4}
\numcell{-2.8}{1.7}{1}
\numcell{-3.3}{4.8}{2}
\numcell{3.3}{4.8}{3}
\numcell{2.8}{1.7}{4}
\cell{-13}{16}{c}{I_2}
\cell{13}{16}{c}{I_2}
\cell{0}{2}{b}{I_3}
\end{picture}
\end{array}
\end{equation}
where $I_3$ is another space with 4 distinguished points
(Figure~\ref{fig:Julia}(c)).  Next, $I_3$ is two copies of itself glued
together:
\begin{equation}
\label{eq:Julia3}
I_3 = 
\begin{array}{c}
\setlength{\unitlength}{1mm} % eps file is 31 x 20 mm
\begin{picture}(32,20)(-16,0)
\cell{0}{0}{b}{\includegraphics{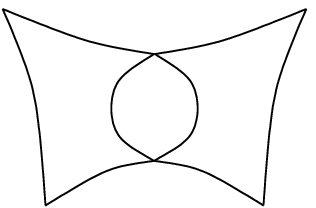}}
\numcell{-3}{5.5}{1}
\numcell{-10}{2.5}{2}
\numcell{-13}{18}{3}
\numcell{-3}{15}{4}
\numcell{3}{5.5}{1}
\numcell{10}{2.5}{2}
\numcell{13}{18}{3}
\numcell{3}{15}{4}
\cell{-8}{10}{c}{I_3}
\cell{8}{10}{c}{I_3}
\end{picture}
\end{array}
\end{equation}
No new spaces appear at this stage, so the process ends.  However, the
one-point space has played a role (since we are gluing at single points),
so I write $I_0$ for the one-point space and record the trivial
equation
\begin{equation}
\label{eq:Julia0}
I_0 = I_0.
\end{equation}

Equations \bref{eq:Julia1}--\bref{eq:Julia0} are a system of simultaneous
equations in which the right-hand sides are `two-dimensional formulas'
expressing each object $I_n$ as a gluing of ($I_m$)s.  Informally, a
self-similarity system is a system of equations like this, and the
particular spaces $I_n$ defined above are a solution of the system.

As a first step towards formalization, observe that the spaces $I_n$
together with their distinguished points define a functor from the category  
\[
\scat{A} = 
\left( \ \ 
\begin{diagram}[height=1.7em,width=3em,noPS]
	&	&1	\\
	&	&	\\
0	&\pile{\rTo\\ \rTo\\ \rTo\\ \rTo}	
		&2	\\
	&\pile{\rdTo\\ \rdTo\\ \rdTo\\ \rdTo}
		&	\\
	&	&3	\\
\end{diagram}
\ \ 
\right)
\]
to the category $\Set$ of sets (or a category of spaces, but let us be
conservative for now).  This describes the left-hand sides of equations
\bref{eq:Julia1}--\bref{eq:Julia0}.  Next, the gluing formulas define a
functor $ G: \ftrcat{\scat{A}}{\Set} \go \ftrcat{\scat{A}}{\Set} $, where
if $X \in \ftrcat{\scat{A}}{\Set}$ then
\[
(GX)_1	
=
\begin{array}{c}
\setlength{\unitlength}{1mm} % eps file is 17 x 20 mm
\begin{picture}(18,20)(-9,-10)
\cell{0.1}{0.1}{c}{\includegraphics{X1new.eps}}
\numcell{-7.7}{2}{1}
\numcell{-0.9}{1.8}{2}
\numcell{0.9}{1.8}{3}
\numcell{7.6}{2}{4}
\numcell{-7.7}{-2}{1}
\numcell{-0.9}{-1.8}{2}
\numcell{0.9}{-1.8}{3}
\numcell{7.6}{-2}{4}
\cell{0}{6}{c}{X_2}
\cell{0}{-6}{c}{X_2}
\end{picture}
\end{array}
=
(X_2 + X_2)/\sim,
\]
and similarly for $(GX)_2$, $(GX)_3$, and $(GX)_0$.  (The picture of
$(GX)_1$ is drawn as if $X_0$ were a single point.)  The simultaneous
equations \bref{eq:Julia1}--\bref{eq:Julia0} assert precisely that
$
I \iso GI
$: 
$I$ is a fixed point of $G$.  

Although these simultaneous equations have many solutions ($G$ has many
fixed points), $I$ is in some sense the universal one.  This means that the
simple diagrams \bref{eq:Julia1}--\bref{eq:Julia0} contain just as much
information as the apparently very complex spaces in
Figure~\ref{fig:Julia}: given the system of equations, we recover these
spaces as the universal solution.  (Caveat: we consider only the intrinsic,
topological aspects of self-similar spaces, not how they are embedded into
an ambient space or any metric structure.)

\begin{figure}
\centering
\setlength{\unitlength}{1mm} % eps file is 51 x 31.5 mm
\begin{picture}(52,32)(-26,-2)
\cell{0.2}{-1}{b}{\includegraphics{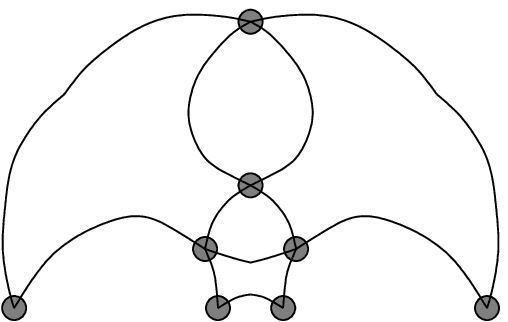}}
\numcell{-23.6}{2.5}{1}
\numcell{-5}{7.8}{2}
\numcell{-2.3}{12.5}{3}
\numcell{-3}{28.5}{4}
\numcell{23.8}{2.5}{1}
\numcell{5.2}{7.8}{2}
\numcell{2.5}{12.5}{3}
\numcell{3}{28.5}{4}
\numcell{-2.8}{1.7}{1}
\numcell{-3.3}{4.8}{2}
\numcell{3.3}{4.8}{3}
\numcell{2.8}{1.7}{4}
\cell{-5.5}{-2}{br}{m}
\cell{-6.5}{5.5}{r}{m'}
\cell{6.5}{5.5}{l}{m''}
\end{picture}
\caption{The eight elements of $M(0, 2)$, with three picked out
(see~\bref{eq:module-functoriality}).} 
\label{fig:module}
\end{figure}

Next we have to make rigorous the notion of `gluing formula', so far
expressed in pictures.  We have a small category $\scat{A}$ whose objects
index the spaces involved, and I claim that the system of gluing formulas
is described by a functor $M: \scat{A}^\op \times \scat{A} \go \Set$.  The
idea is that for $b, a \in \scat{A}$,
\begin{eqnarray*}
M(b, a)		&=	&
\{\textrm{copies of the }b\textrm{th space used in the gluing formula} \\
		&	&
\ \,\textrm{for the }a\textrm{th space} \}.
\end{eqnarray*}
For example, in the gluing formula for $I_2$, the one-point space $I_0$
appears 8 times (Figure~\ref{fig:module}),
$I_1$ does not appear at all, $I_2$ appears twice, and $I_3$ appears once,
so, writing $n$ for an $n$-element set,
\[
M(0, 2) = 8,
\diagspace
M(1, 2) = \emptyset,
\diagspace
M(2, 2) = 2,
\diagspace
M(3, 2) = 1.
\]
Similarly, the gluing formula for $I_0$ is nothing but a single copy of
$I_0$, so
\[
M(0, 0) = 1,
\diagspace
M(1, 0) = M(2, 0) = M(3, 0) = \emptyset.
\]
For functoriality, consider, for instance, the second of the four maps $0
\go 2$.  (In the notation above, this determines the second of the four
distinguished points of $I_2$.)  This induces functions
\begin{equation}
\label{eq:module-functoriality}
M(0, 0) \go M(0, 2),
\diagspace
M(2, 2) \go M(0, 2)
\end{equation}
(among others).  The first function sends the unique element of $M(0, 0)$
to the element of $M(0, 2)$ marked $m$ in Figure~\ref{fig:module}.  Writing
$M(2, 2) = \{ \textrm{Left}, \textrm{ Right} \}$, the second function sends
$\textrm{Left}$ to $m'$ and $\textrm{Right}$ to $m''$.  

It is convenient to use the language of modules.  Given small categories
$\scat{A}$ and $\scat{B}$, a \demph{module} $M: \scat{B} \gomod \scat{A}$
is a functor $M: \scat{B}^\op \times \scat{A} \go \Set$.  For example,
suppose that $\scat{A}$ and $\scat{B}$ are monoids ($=$ categories with
only one object): then a module $\scat{B} \gomod \scat{A}$ is a set with a
compatible left $\scat{A}$-action and right $\scat{B}$-action.  (If we work
with categories enriched in abelian groups then one-object categories are
rings and modules are bimodules.)  Write $b \gobymod{m} a$ to mean $m \in
M(b, a)$: then a module $M: \scat{B} \gomod \scat{A}$ is an indexed family
$(M(b, a))_{b \in \scat{B}, a \in \scat{A}}$ of sets together with actions
\begin{eqnarray*}
b \gobymod{m} a \goby{f} a'	&
\ \goesto	&\ 
b \gobymod{fm} a',	\\
b' \goby{g} b \gobymod{m} a	&
\ \goesto	&\ 
b' \gobymod{mg} a	
\end{eqnarray*}
such that $(f'f)m = f'(fm)$, $1m = m$, and dually, and $(fm)g = f(mg)$.

(Experts should note that I have reversed the standard convention: a
functor $M: \scat{B}^\op \times \scat{A} \go \Set$ is written $\scat{B}
\gomod \scat{A}$, not $\scat{A} \gomod \scat{B}$.  One reason is that given
such an $M$, an element $m \in M(b, a)$ can sensibly be written as $b
\gobymod{m} a$; then with the convention used here, a module $\scat{B}
\gobymod{M} \scat{A}$ consists of elements $b \gobymod{m} a$.)

The system of self-similarity equations \bref{eq:Julia1}--\bref{eq:Julia0}
is therefore encoded by a small category $\scat{A}$ and a module $M:
\scat{A} \gomod \scat{A}$.  The conceptual distinction between arrows $b
\go a$ in $\scat{A}$ and arrows $b \gomod a$ in $M$ is that the arrows in
$\scat{A}$ determine where gluing may potentially take place, but the
arrows in $\scat{M}$ say what the gluing actually is.

Here is another way of seeing that the self-similarity equations are
encoded by a module.  The right-hand side of each of
\bref{eq:Julia1}--\bref{eq:Julia0} is a formal gluing of objects of
$\scat{A}$.  `Gluings' are colimits, so if $\widehat{\scat{A}}$ is the
category obtained by taking $\scat{A}$ and freely adjoining all possible
colimits then the system of equations amounts to a functor from $\scat{A}$
to $\widehat{\scat{A}}$.  But $\widehat{\scat{A}} = \pshf{\scat{A}}$
(see~\cite[I.5.4]{MM}), so the system is a functor $\scat{A} \go
\pshf{\scat{A}}$, that is, a module $\scat{A} \gomod \scat{A}$.

Simple colimits can often be described by diagrams or formulas, which
provides a useful way of specifying simple self-similarity systems.  For
instance, equations \bref{eq:Julia1}--\bref{eq:Julia0} are an informal
description of our $(\scat{A}, M)$.  In the same way, the self-similarity
system $(\scat{A}, M)$ in which $\scat{A}$ is the discrete category with
object-set $\{ 0, 1 \}$ and $M(0, 0) = M(0, 1) = M(1, 1) = 1$ and $M(1, 0) =
\emptyset$ can informally (and more intelligibly) be described by the
equations
\begin{eqnarray}
A	&=	&A		\label{eq:informal-SSS}	\\
B	&=	&A + B.		\nonumber
\end{eqnarray}

Given rings $A$, $B$ and $C$, an $(A, B)$-bimodule $M$, and a $(B,
C)$-bimodule $N$, there arises an $(A, C)$-bimodule $M \otimes_B N$.  There
is a similar tensor product of categorical modules:
\[
\scat{C} \gobymod{N} \scat{B} \gobymod{M} \scat{A}
\diagspace
\textrm{gives rise to}
\diagspace
\scat{C} \gobymod{M \otimes N} \scat{A}.
\]
Here $M \otimes N$ is defined by the coend formula
\[
(M \otimes N) (c, a)
=
\int^b M(b, a) \times N(c, b).
\]
(See \cite[Ch.~IX]{CWM} for an explanation of coends.)  Concretely,
\[
(M \otimes N) (c, a)
=
\left(
\sum_b M(b, a) \times N(c, b)
\right)
/\sim
\]
where $\sim$ is the equivalence relation generated by $(mg, n) \sim (m,
gn)$ for all $m$, $g$, $n$ of the appropriate types; the element of $(M
\otimes N) (c, a)$ represented by $(m, n) \in M(b, a) \times N(c, b)$ is
written $m \otimes n$.  The tensor product of modules is associative and
unital up to coherent isomorphism.  (More precisely, categories, modules,
and their maps form a bicategory: \cite[7.8.2]{Bor1}.)

In the Julia set example, the endofunctor $G$ of $\ftrcat{\scat{A}}{\Set}$
is just $M \otimes \dashbk$.  This makes sense: a functor $X: \scat{A} \go
\Set$ can be regarded as a module $X: \One \gomod \scat{A}$ (where $\One$
is the category with one object and only the identity arrow), and there is
then a tensor product $M \otimes X: \One \gomod \scat{A}$, that is, a
functor $M \otimes X: \scat{A} \go \Set$.  It is given by
\begin{equation}
\label{eq:mod-act}	
(M \otimes X) a	
=
\int^b M(b, a) \times Xb	
=	
\left(
\sum_b M(b, a) \times Xb
\right)
/\sim.		
\end{equation}
In the example, taking $a = 2$, this says that
\[
(M \otimes X)_2
=
( 8 \times X_0 + 2 \times X_2 + X_3 )
/\sim
\]
where $\sim$ identifies the various copies of $X_0$ with their images in
$X_2$ and $X_3$.  In general,~\bref{eq:mod-act} makes precise the idea that
$M(b, a)$ is the set of copies of the $b$th space used in the gluing
formula for the $a$th space.

The second example is a result of Peter Freyd, and comes from a very
different direction.  To state it we need some more terminology.

Given a category $\cat{C}$ and an endofunctor $G$ of $\cat{C}$, a
\demph{$G$-coalgebra} is an object $X$ of $\cat{C}$ together with a map
$\xi: X \go GX$.  (For instance, if $\cat{C}$ is a category of modules and
$GX = X \otimes X$ then a $G$-coalgebra is a coalgebra---not necessarily
coassociative---in the usual sense.)  A \demph{map} $(X, \xi) \go (X',
\xi')$ of coalgebras is a map $X \go X'$ in $\cat{C}$ making the evident
square commute.  Depending on what $G$ is, the category of $G$-coalgebras
may or may not have a terminal object, but if it does then it is a fixed
point:
\begin{lemma}[Lambek~\cite{Lam}]
\label{lemma:Lambek}
Let $\cat{C}$ be a category and $G$ an endofunctor of $\cat{C}$.  If $(I,
\iota)$ is terminal in the category of $G$-coalgebras then $\iota: I \go
GI$ is an isomorphism.
\done
\end{lemma}

Here is what Freyd said, strengthened slightly.  Let $\cat{C}$ be the category
whose objects are diagrams $X_0 \parpair{u}{v} X_1$ where $X_0$ and $X_1$
are sets and $u$ and $v$ are injections with disjoint images; then an
object of $\cat{C}$ can be drawn as
\[
\setlength{\unitlength}{1mm}
\begin{picture}(30,14)(-15,-5)
\cell{0}{0}{c}{\includegraphics{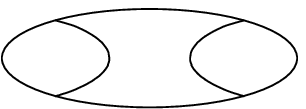}}
\cell{-10}{0}{c}{X_0}
\cell{10}{0}{c}{X_0}
\cell{0}{5.5}{b}{X_1}
\end{picture}
\]
where the copies of $X_0$ on the left and the right are the
images of $u$ and $v$ respectively.  A map $X \go X'$ in
$\cat{C}$ consists of functions $X_0 \go X'_0$ and $X_1 \go X'_1$
making the evident two squares commute.  Now, given $X \in
\cat{C}$ we can form a new object $GX$ of $\cat{C}$ by gluing
two copies of $X$ end to end:
\begin{equation}
\label{eq:Freyd-gluing}
\begin{array}{c}
\setlength{\unitlength}{1mm}
\begin{picture}(49,14)(-24.5,-5)
\cell{0}{0}{c}{\includegraphics{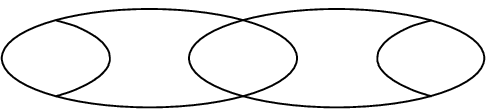}}
\cell{0}{0}{c}{X_0}
\cell{-20}{0}{c}{X_0}
\cell{20}{0}{c}{X_0}
\cell{-10}{5.5}{b}{X_1}
\cell{10}{5.5}{b}{X_1}
\end{picture}
\end{array}
.
\end{equation}
Formally, the endofunctor $G$ of $\cat{C}$ is defined by pushout:
\[
\begin{diagdiag}
	&	&	&	&(GX)_1	&	&	&	&	\\
	&	&	&\ruTo	&	&\luTo	&	&	&	\\
	&	&X_1	&	&\textrm{pushout} &
						&X_1	&	&	\\
	&\ruTo<u&	&\luTo>v&	&\ruTo<u&	&\luTo>v&	\\
(GX)_0 = X_0&	&	&	&X_0	&	&	&	&X_0.	\\
\end{diagdiag}
\]
For example, the unit interval with its endpoints distinguished forms an
object 
\[
I = \left( \{\star\} \parpair{0}{1} [0,1] \right)
\]
of $\cat{C}$, and $GI$ is naturally described as an interval of
length 2, again with its endpoints distinguished:
\[
GI = \left( \{\star\} \parpair{0}{2} [0,2] \right).
\]
So there is a coalgebra structure $\iota: I \go GI$ on $I$ given by
multiplication by two.  Freyd's Theorem says that this is, in fact, the
\emph{universal} example of a $G$-coalgebra:
\begin{thm}[Freyd~\cite{Fre}]	\label{thm:Freyd}
$(I, \iota)$ is terminal in the category of $G$-coalgebras.
\end{thm}

This follows from a general result~\cite[\egrecFreyd]{SS2}.  A
direct proof runs roughly as follows.  Take a $G$-coalgebra $(X,
\xi)$ and an element $x_0 \in X_1$.  Then $\xi(x_0) \in (GX)_1$
is in either the left-hand or the right-hand copy of $X_1$, so
gives rise to a binary digit $m_1 \in \{0, 1\}$ and a new element
$x_1 \in X_1$.  (If $\xi(x_0)$ is in the overlap between the two
copies of $X_1$, choose left or right arbitrarily.)  Iterating
gives a binary representation $0.m_1 m_2 \ldots$ of an element of
$[0, 1]$, and this is the image of $x_0$ under the unique
coalgebra map $(X, \xi) \go (I, \iota)$.

The formalism introduced in the Julia set example can also be used here.
The category $\cat{C}$ is a full subcategory of $\ftrcat{\scat{A}}{\Set}$,
where $\scat{A} = (0 \parpair{\sigma}{\tau} 1)$.  The module $M: \scat{A}
\gomod \scat{A}$ corresponding to the gluing formula~\bref{eq:Freyd-gluing}
is given by
\begin{equation}	\label{eq:Freyd-SSS}
\begin{diagram}
				&\ \ \				&
M(\dashbk, 0)			&\pile{\rTo^{\sigma\cdot\dashbk}\\ 
				\rTo_{\tau\cdot\dashbk}}	&
M(\dashbk, 1)			\\
				&\				&
				&				&
				\\
M(0, \dashbk)			&				&
\{ \id \}			&\pile{\rTo^0\\ \rTo_1}		&
\{ 0, \dhalf, 1 \}		\\
\uTo<{\dashbk\cdot\sigma} 
\uTo>{\dashbk\cdot\tau}		&				&
\uTo \uTo			&				&
\uTo<\inf \uTo>\sup		\\
M(1, \dashbk)			&				&
\emptyset			&\pile{\rTo\\ \rTo}		&
\{ [0, \dhalf], [\dhalf, 1] \}.	\\
\end{diagram}
\end{equation}
(Here $M(0, 1)$ is just a $3$-element set and $M(1, 1)$ a $2$-element set,
but their elements have been named suggestively.)  Then $M \otimes \dashbk$
defines an endofunctor of $\ftrcat{\scat{A}}{\Set}$, and the endofunctor
$G$ of $\cat{C}$ is its restriction.

The only remaining mystery is the condition that the functions $u, v: X_0
\go X_1$ are injections with disjoint images.  Without it, the theorem
would degenerate entirely, as the terminal coalgebra would be $(\{ \star \}
\parpairu \{ \star \})$.  It turns out to be a form of flatness.

First note that any two functors $X: \scat{A} \go \Set$ and $Y:
\scat{A}^\op \go \Set$ on a small category $\scat{A}$ have a tensor product
$Y \otimes X$ (a mere set), since they can be regarded as modules
\[
\One \gobymod{X} \scat{A} \gobymod{Y} \One.
\]
Explicitly,
\[
Y \otimes X
=
\int^a Ya \times Xa
=
\left(
\sum_a Ya \times Xa
\right) 
/\sim.
\]

By definition, a left module $X$ over a ring is flat if the functor
$\dashbk \otimes X$ preserves finite limits.  There is an analogous
definition when $X$ is a $\Set$-valued functor, but we actually want
something weaker:
\begin{defn}
Let $\scat{A}$ be a small category.  A functor $X: \scat{A} \go \Set$ is
\demph{componentwise flat} or \demph{nondegenerate} if the functor
\[
\dashbk \otimes X: \pshf{\scat{A}} \go \Set
\]
preserves finite connected limits.  The full subcategory of
$\ftrcat{\scat{A}}{\Set}$ consisting of the nondegenerate functors is
written $\ndSet{\scat{A}}$.
\end{defn}
(A category is \demph{connected} if it is nonempty and cannot be written as
a coproduct of two nonempty categories; a \demph{finite connected limit} is
a limit over a finite connected category.)

This definition is explained at length in~\S\ref{sec:nondegen}.  It is
shown there that a $\Set$-valued functor $X$ on $(0 \parpair{\sigma}{\tau}
1)$ is nondegenerate exactly when the functions $X\sigma$ and $X\tau$ are
injections with disjoint images: so $\cat{C} = \ndSet{\scat{A}}$.
Proposition~\ref{propn:preservationofnondegeneracy} says that for general
$\scat{A}$ and $M$, the endofunctor $M \otimes \dashbk$ of
$\ftrcat{\scat{A}}{\Set}$ restricts to an endofunctor of $\ndSet{\scat{A}}$
as long as $M$ is nondegenerate:
\begin{defn}
Let $\scat{A}$ and $\scat{B}$ be small categories.  A module $M: \scat{B}
\gomod \scat{A}$ is \demph{nondegenerate} if $M(b, \dashbk): \scat{A} \go
\Set$ is nondegenerate for each $b \in \scat{B}$. 
\end{defn}

It is visible from diagram~\bref{eq:Freyd-SSS} that the Freyd module $M$ is
nondegenerate.

To summarize: starting with a certain small category $\scat{A}$ and a
certain nondegenerate module $M: \scat{A} \gomod \scat{A}$, Freyd's Theorem
describes the terminal coalgebra for the endofunctor $M \otimes \dashbk$ of
$\ndSet{\scat{A}}$.  The Julia set example uses a different $\scat{A}$ and
$M$.  In general, I will restrict to those $\scat{A}$ and $M$ for which
`each gluing formula is finite', although I allow there to be infinitely
many such formulas (infinitely many objects of $\scat{A}$).

To make this precise, recall that any presheaf $Y: \scat{A}^\op \go \Set$
on a small category $\scat{A}$ has a \demph{category of elements}%
\label{p:cat-elts}
$\elt{Y}$, whose objects are pairs $(a, y)$ with $a \in \scat{A}$ and $y
\in Ya$ and whose maps $(a, y) \go (a', y')$ are maps $f: a \go a'$ in
$\scat{A}$ such that $(Yf)y' = y$.  Similarly, any covariant functor $X:
\scat{A} \go \Set$ has a \demph{category of elements} $\elt{X}$.  In each
case, there is a covariant projection functor from the category of elements
to $\scat{A}$.

A module $M: \scat{A} \rMod \scat{A}$ is \demph{finite} if for each $a \in
\scat{A}$, the category $\elt{M(\dashbk, a)}$ is finite.  Explicitly, this
says that for each $a \in \scat{A}$ there are only finitely many diagrams
of the form
\[
b' \goby{f} b \gobymod{m} a.
\]
Certainly this holds if, as in the Freyd example, the category $\scat{A}$
and the sets $M(b, a)$ are finite.

\begin{defn}
A \demph{self-similarity system} is a small category $\scat{A}$ together
with a finite nondegenerate module $M: \scat{A} \gomod \scat{A}$.
\end{defn}

\begin{defn}
Let $(\scat{A}, M)$ be a self-similarity system.  An \demph{$M$-coalgebra
(in $\Set$)} is a coalgebra for the endofunctor $M \otimes \dashbk$ of
$\ndSet{\scat{A}}$.  A \demph{universal solution} of $(\scat{A}, M)$ (in
$\Set$) is a terminal object in the category of $M$-coalgebras.
\end{defn}

In this language, Freyd's Theorem describes the universal solution of a
certain self-similarity system.

There is also a topological version of Freyd's Theorem.  One's first
thought might be to take the definitions of $\cat{C}$ and $G$ and change
`set' to `space' and `function' to `continuous map'; but then the universal
solution is given by the indiscrete topology on $[0, 1]$.  The Euclidean
topology appears, however, if we insist that $u, v: X_0 \go X_1$ are closed
maps.  So, let $\cat{C}'$ be the category whose objects are diagrams $X_0
\parpairu X_1$ of topological spaces and continuous closed injections with
disjoint images, define an endofunctor $G'$ of $\cat{C}'$ just as $G$ was
defined, and define $(I, \iota)$ as before, with the Euclidean topology on
$[0, 1]$.  Then:
\begin{thm}[Topological Freyd]
\label{thm:topologicalFreyd}
$(I, \iota)$ is terminal in the category of $G'$-coalgebras.
\end{thm}
(See~\cite[\egrecFreyd]{SS2} for a proof.)  This suggests the following
general definition.
\begin{defn}
Let $\scat{A}$ be a small category.  Write $U: \Top \go \Set$ for the
underlying set functor.  A functor $X: \scat{A} \go \Top$ is
\demph{nondegenerate} if $U \of X$ is nondegenerate and for each map $f$ in
$\scat{A}$, the map $Xf$ is closed.  The full subcategory of
$\ftrcat{\scat{A}}{\Top}$ formed by the nondegenerate functors is written
$\ndTop{\scat{A}}$. 
\end{defn}
We will often meet functors $X$ for which each space $Xa$ is compact
Hausdorff, and then the closedness condition is automatically satisfied.

We also want a general notion of topological $M$-coalgebra.  Let $\cat{E}$
be a category with finite colimits, $\scat{A}$ a small category, and $M:
\scat{A} \gomod \scat{A}$ a finite module.  Then there is an endofunctor $M
\otimes \dashbk$ of $\ftrcat{\scat{A}}{\cat{E}}$ defined by the usual coend
formula
\[
(M \otimes X) a 
=
\int^b M(b, a) \times Xb.
\]
Proposition~\ref{propn:topologicalmoduleaction} says that for any
self-similarity system $(\scat{A}, M)$, the endofunctor $M \otimes \dashbk$
of $\ftrcat{\scat{A}}{\Top}$ restricts to an endofunctor of
$\ndTop{\scat{A}}$.

\begin{defn}
Let $(\scat{A}, M)$ be a self-similarity system.  An \demph{$M$-coalgebra
in $\Top$} is a coalgebra for the endofunctor $M \otimes \dashbk$ of
$\ndTop{\scat{A}}$.  A \demph{universal solution of $(\scat{A}, M)$ in
$\Top$} is a terminal $M$-coalgebra in $\Top$.
\end{defn}

For example, the topological Freyd theorem describes the universal solution
in $\Top$ of a certain self-similarity system. 

Universal solutions are evidently unique (up to canonical isomorphism) when
they exist.  The word `solution' is justified by Lambek's
Lemma~(\ref{lemma:Lambek}): if $(J, \gamma)$ is a universal solution then
$M \otimes J \iso J$.  Note that the converse fails: the empty functor $J =
\emptyset$ has a unique coalgebra structure, is nondegenerate, and
satisfies $M \otimes J \iso J$, but is not usually the terminal coalgebra.
(It is the initial algebra.)  In the Freyd interval example, there are many
coalgebras $(J, \gamma)$ such that $\gamma$ is an isomorphism but $(J,
\gamma)$ is not the universal solution: for instance, 
the universal
solution can be multiplied by any space $S$ to give
such a coalgebra $(J, \gamma)$, with $J = \left(S \parpairu [0, 1]
\times S\right)$.

Just as an ordinary system of equations need not have a solution, a
self-similarity system need not have a universal solution.
In~\S\ref{sec:univ-soln} we meet an explicit condition equivalent to the
existence of a universal solution.

\section{Nondegeneracy}
\label{sec:nondegen}

In this section I explain nondegeneracy, first by theory and then by
examples.  The theory leads up to the result that a functor $X: \scat{A}
\go \Set$ is nondegenerate if and only if it satisfies the following
explicit conditions:
\begin{description}
\item[\cstyle{ND1}]
given
\[
\begin{diagdiag}
a	&		&	&		&a'	\\
	&\rdTo<f	&	&\ldTo>{f'}	&	\\
	&		&b	&		&	\\
\end{diagdiag}
\]
in $\scat{A}$ and $x \in Xa$, $x' \in Xa'$ such that $fx = f'x'$, there
exist a commutative square
\[
\begin{diagdiag}
	&		&c	&		&	\\
	&\ldTo<g	&	&\rdTo>{g'}	&	\\
a	&		&	&		&a'	\\
	&\rdTo<f	&	&\ldTo>{f'}	&	\\
	&		&b	&		&	\\
\end{diagdiag}
\]
and $z \in Xc$ such that $x = gz$, $x' = g'z$
\item{\cstyle{ND2}} given $a \parpair{f}{f'} b$ in $\scat{A}$ and $x \in
Xa$ such that $fx = f'x$, there exist a fork 
\begin{equation}	\label{eq:ND2-fork}
c \goby{g} a \parpair{f}{f'} b
\end{equation}
and $z \in Xc$ such that $x = gz$.  
(A diagram~\bref{eq:ND2-fork} is a \demph{fork} if $fg = f'g$.)
\end{description}
The examples illustrate that nondegeneracy means `no unforced equalities'.

For the theory I will assume some more sophisticated categorical knowledge
than in the rest of the paper.  Readers who prefer to take it on trust can
jump to the passage after
Corollary~\ref{cor:categorywithpullbacksandequalizers}.

None of this theory is new: it goes back to Grothendieck and
Verdier~\cite{GV} and Gabriel and Ulmer~\cite{GU}, and was later developed
by Weberpals~\cite{Web}, Lair~\cite{Lai}, Ageron~\cite{Age}, and Ad\'amek,
Borceux, Lack, and Rosick\'y~\cite{ABLR}.  More general statements of much
of what follows can be found in~\cite{ABLR}.

Let us begin with ordinary flat functors.  A functor $X: \scat{A} \go \Set$
on a small category $\scat{A}$ is \demph{flat} if $\dashbk \otimes X:
\pshf{\scat{A}} \go \Set$ preserves finite limits.

\begin{thm}[Flatness]
\label{thm:flatness}
Let $\scat{A}$ be a small category.  The following conditions on a functor
$X: \scat{A} \go \Set$ are equivalent:
\begin{enumerate}
\item \label{item:flat-flat}
$X$ is flat
\item \label{item:flat-fin}
every finite diagram in $\elt{X}$ admits a cone
\item \label{item:flat-fin-bits}
each of the following holds:
\begin{itemize}
\item there exists $a \in \scat{A}$ for which $Xa$ is nonempty
\item given $a, a' \in \scat{A}$, $x \in Xa$, and $x' \in
  Xa'$, there exist a diagram $a \ogby{g} c \goby{g'} a'$ in $\scat{A}$ and
  $z \in Xc$ such that $gz = x$ and $g'z = x'$
\item \cstyle{ND2}.
\end{itemize}
\end{enumerate}
\end{thm}
\begin{proof}
See \cite[\S 6.3]{Bor1} or \cite[VII.6]{MM}, for instance.  
\done
\end{proof}

The following lemmas are often used to prove this theorem and will also be
needed later.
\begin{lemma}[Existence of cones]
\label{lemma:existenceofcones}
Let $\scat{I}$ and $\scat{A}$ be small categories and let $X: \scat{A} \go
\Set$.  If $\dashbk \otimes X: \pshf{\scat{A}} \go \Set$ preserves limits
of shape $\scat{I}$ then every diagram of shape $\scat{I}$ in $\elt{X}$
admits a cone.
\end{lemma}
\paragraph*{Remark} The hypothesis can be weakened to say that $\dashbk
\otimes X$ preserves limits of shape $\scat{I}$ \emph{of diagrams of
representables}, that is, of diagrams $\scat{I} \go \pshf{\scat{A}}$ that
factor through the Yoneda embedding of $\scat{A}$.

\begin{proof}
Let $D: \scat{I} \go \elt{X}$ be a diagram of shape $\scat{I}$, writing
$Di = (a_i, x_i)$ for each $i \in \scat{I}$.  Then there is a diagram
$Y_\blob: \scat{I} \go \pshf{\scat{A}}$ given by $Y_i = \scat{A}(\dashbk,
a_i)$, so by hypothesis the canonical map
\[
\int^a \left( \int_i \scat{A}(a, a_i) \right) \times Xa 
\iso
\left( \int_i Y_i \right) \otimes X
\go
\int_i (Y_i \otimes X)
\iso
\int_i Xa_i
\]
is a bijection, and in particular a surjection.  Since $(x_i)_{i \in
\scat{I}} \in \int_i Xa_i$, there exist $a \in \scat{A}$ and
\[
((p_i)_{i \in \scat{I}}, x) 
\in 
\left( \int_i \scat{A}(a, a_i) \right) \times Xa
\]
such that $p_i x = x_i$ for all $i$.  Hence 
$
\left( 
(a, x) \goby{p_i} (a_i, x_i)
\right)_{i \in \scat{I}}
$
is a cone on $D$.
\done
\end{proof}

Say that a category $\cat{C}$ has the \demph{square-completion property} if
there exists a cone on every diagram of shape $\littlepullback$ in
$\cat{C}$.

\begin{lemma}[Connectedness by spans]
\label{lemma:connectednessbyspans}
Two objects $c, c'$ of a category with the square-completion property are
in the same connected-component if and only if there exists a span
$
c \og c'' \go c'
$
connecting them.
\done
\end{lemma}

\begin{lemma}[Equality in a tensor product]
\label{lemma:equalityinatensorproduct}
Let $\scat{A}$ be a small category and 
\[
X: \scat{A} \go \Set,
\diagspace
Y: \scat{A}^\op \go \Set.
\]
Suppose that $\elt{X}$ has the square-completion property.  Let 
\[
a, a' \in \scat{A},
\diagspace
(y, x) \in Ya \times Xa,
\diagspace
(y', x') \in Ya' \times Xa'.
\]
Then $y \otimes x = y' \otimes x'$ if and only if there exists a span
$
a \ogby{f} b \goby{f'} a'
$
and an element $z \in Xb$ such that $x = fz$, $x' = f'z$, and $yf = y'f'$.
\end{lemma}
\begin{proof}
See the remarks after the statement of Theorem VII.6.3 in~\cite{MM}.
\done
\end{proof}

There is a characterization of componentwise flat functors very similar to
that of flat functors in Theorem~\ref{thm:flatness}.  First we need a
fact about connectedness.

\begin{lemma}[Components of a functor]
\label{lemma:componentsofafunctor}
Any functor $X: \scat{A} \go \Set$ on a small category $\scat{A}$ can be
written as a sum $X \iso \sum_{j \in J} X_j$ where $J$ is some set and
$\elt{X_j}$ is connected for each $j\in J$.
\end{lemma}

\begin{proof}
We use the equivalence between $\Set$-valued functors and discrete
opfibrations.  Write $\elt{X}$ as a sum $\sum_{j \in J} \scat{E}_j$ of
connected categories.  For each $j$, the restriction to $\scat{E}_j$ of the
projection $\elt{X} \go \scat{A}$ is still a discrete opfibration, so
corresponds to a functor $X_j: \scat{A} \go \Set$.  Then
\[
\elt{\sum X_j}
\iso 
\sum \elt{X_j}
\iso
\sum \scat{E}_j
\iso
\elt{X}
\]
compatibly with the projections, so $\sum X_j \iso X$.  \done
\end{proof}

Here is the main result.
\begin{thm}[Componentwise flatness]
\label{thm:componentwiseflatness}
Let $\scat{A}$ be a small category.  The following conditions on a functor
$X: \scat{A} \go \Set$ are equivalent:
\begin{enumerate}
\item \label{item:cwflat-cwflat}
$X$ is componentwise flat
\item \label{item:cwflat-finconn}
every finite connected diagram in $\elt{X}$ admits a cone
\item \label{item:cwflat-finconn-bits}
$X$ satisfies \cstyle{ND1} and~\cstyle{ND2}
\item \label{item:cwflat-sum}
$X$ is a sum of flat functors.
\end{enumerate}
\end{thm}
\paragraph*{Remark} In Lemma~\ref{lemma:componentsofafunctor}, the functors
$X_j$ may be regarded as the connected-components of $X$.  A further
equivalent condition is that every connected-component of $X$ is flat:
hence the name `componentwise flat'.
\begin{proof}

\paragraph*{\bref{item:cwflat-cwflat}$\implies$\bref{item:cwflat-finconn}}
Follows from Lemma~\ref{lemma:existenceofcones}.

\paragraph*{\bref{item:cwflat-finconn}$\implies$\bref{item:cwflat-finconn-bits}}
\cstyle{ND1} says that every diagram of shape $\littlepullback$ in
$\elt{X}$ admits a cone, and similarly~\cstyle{ND2} for $\littleequalizer$.

\paragraph*{\bref{item:cwflat-finconn-bits}$\implies$\bref{item:cwflat-sum}}
Write $X \iso \sum_{j\in J} X_j$ as in
Lemma~\ref{lemma:componentsofafunctor}.  Then in each $\elt{X_j}$, there
exists a cone on every diagram of shape
\[
\littlepullback
\diagspace
\textrm{or}
\diagspace
\littleequalizer
\]
(since $\elt{X} \iso \sum_j \elt{X_j}$), of shape $\emptyset$ (since
$\elt{X_j}$ is connected and therefore nonempty), and of shape
$\littleproduct$ (since $\elt{X_j}$ is connected and has the
square-completion property).  So by
\bref{item:flat-fin-bits}$\implies$\bref{item:flat-flat} of
Theorem~\ref{thm:flatness}, each $X_j$ is flat.

\paragraph*{\bref{item:cwflat-sum}$\implies$\bref{item:cwflat-cwflat}}
Any sum of componentwise flat functors is componentwise flat, as follows
from the fact that sums commute with connected limits in $\Set$.
\done
\end{proof}

\begin{cor}[Componentwise filtered categories]
\label{cor:componentwisefilteredcategories}
The following conditions on a small category $\scat{B}$ are equivalent:
\begin{enumerate}
\item \label{item:cwfilt-comm}
finite connected limits commute with colimits of shape $\scat{B}$ in $\Set$
\item \label{item:cwfilt-finconn}
every finite connected diagram in $\scat{B}$ admits a cocone
\item \label{item:cwfilt-finconn-bits} 
every diagram $b_1 \og b_3 \go b_2$ in $\scat{B}$ can be completed to a
commutative square, and every parallel pair $b_1 \parpair{f}{f'} b_2$
of arrows in $\scat{B}$ can be extended to a cofork.
\end{enumerate}
\end{cor}
\begin{proof}
In Theorem~\ref{thm:componentwiseflatness}, take $\scat{A} = \scat{B}^\op$
and $X$ to be the functor with constant value $1$.  Then $\elt{X} \iso
\scat{B}^\op$ and $\dashbk \otimes X$ forms colimits.  The result follows.
\done
\end{proof}

A small category $\scat{B}$ satisfying the equivalent conditions of
Corollary~\ref{cor:componentwisefilteredcategories} is called
\demph{componentwise filtered}.  (Grothendieck and Verdier say
`pseudo-filtrante' \cite{GV}, and a further equivalent condition is that
each connected-component is filtered.)  So $X: \scat{A} \go \Set$ is
componentwise flat just when $\elt{X}$ is componentwise \emph{co}filtered.

Componentwise flatness relates to limit preservation as follows.

\begin{lemma}
\label{lemma:preservationbyacomponentwiseflatfunctor}
Componentwise flat functors preserve finite connected limits.
\end{lemma}
\begin{proof}
Let $\scat{A}$ be a small category and $X: \scat{A} \go \Set$ a
componentwise flat functor.  We have
\begin{equation}	\label{eq:Yonedatensor}
X 
\iso
\left(
\scat{A} \goby{\mr{Yoneda}} 
\pshf{\scat{A}} \goby{\dashbk \otimes X} 
\Set
\right)
\end{equation}
and the Yoneda embedding preserves limits.
\done
\end{proof}

\begin{cor}
\label{cor:categorywithpullbacksandequalizers}
Let $\scat{A}$ be a small category with all pullbacks and equalizers.  Then
a functor $X: \scat{A} \go \Set$ is componentwise flat if and only if it
preserves pullbacks and equalizers.
\end{cor}
\begin{proof}
Suppose that $X$ preserves pullbacks and equalizers.
By~\bref{eq:Yonedatensor}, $\dashbk \otimes X$ preserves pullbacks and
equalizers of representables, so by the Remark after the statement of
Lemma~\ref{lemma:existenceofcones}, any diagram of shape $\littlepullback$
or $\littleequalizer$ in $\elt{X}$ admits a cone.  But this says that $X$
satisfies \cstyle{ND1} and \cstyle{ND2}.  \done
\end{proof}

Intuitively, nondegeneracy (componentwise flatness) of a functor $X:
\scat{A} \go \Set$ says that no equation between elements of $X$ holds
unless it must.  For example, if $\scat{A} = (\blob \go \blob)$ then a
functor $X: \scat{A} \go \Set$ is a map $(X_0 \goby{i} X_1)$ of sets, and
nondegeneracy of $X$ says that the equation $ix_0 = ix'_0$ holds only when
it must, that is, only when $x_0 = x'_0$; thus, $X$ is nondegenerate just
when $i$ is injective (Example~\ref{eg:nondegen-arrowcat}).  Or, let
$\scat{A}$ be the category generated by objects and arrows
\begin{equation}
\label{eq:cofork}
0 \parpair{\sigma}{\tau} 1 \goby{\rho} 2
\end{equation}
subject to $\rho\sigma = \rho\tau$.  If $X: \scat{A} \go \Set$ is
nondegenerate then the equation $(X\rho) x_1 = (X\rho) x'_1$ holds only
when it must, that is, when $x_1 = x'_1$ or there exists $x_0$ satisfying
$\{ x_1, x'_1 \} = \{ (X\sigma) x_0, (X\tau) x_0 \}$
(Example~\ref{eg:nondegen-cofork}).

Let us work out what nondegeneracy says for various specific categories
$\scat{A}$.  Note that \cstyle{ND1} holds automatically if either $f$ or
$f'$ is an isomorphism, and that \cstyle{ND2} holds automatically if $f =
f'$; we therefore ignore these cases.  Moreover, if $f$ is monic then
\cstyle{ND1} in the case $f = f'$ just says that $Xf$ is injective;
indeed, we already know from
Lemma~\ref{lemma:preservationbyacomponentwiseflatfunctor} that
nondegenerate functors preserve monics.  

\begin{example} \label{eg:nondegen-arrowcat}
Let $\scat{A} = \left( 0 \goby{\sigma} 1 \right)$.  Then $X: \scat{A} \go
\Set$ is nondegenerate if and only if the function $X\sigma: X0 \go X1$ is
injective. 
\end{example}

\begin{example} \label{eg:nondegen-Freyd}
Let $\scat{A} = \left( 0 \parpair{\sigma}{\tau} 1 \right)$, so that a
functor $X: \scat{A} \go \Set$ is a pair $\left( X_0
\parpair{X\sigma}{X\tau} X_1 \right)$ of functions.  Then \cstyle{ND1} in
the case $f = f'$ says that $X\sigma$ and $X\tau$ are injective.  The only
other nontrivial case of \cstyle{ND1} is $f = \sigma$, $f' = \tau$, and
since the diagram
\[
\begin{diagdiag}
0	&		&	&		&0	\\
	&\rdTo<\sigma	&	&\ldTo>\tau	&	\\
	&		&1	&		&	\\
\end{diagdiag}
\]
cannot be completed to a commutative square, \cstyle{ND1} says that
$X\sigma$ and $X\tau$ have disjoint images.  The only nontrivial case of
\cstyle{ND2} is $f = \sigma$, $f' = \tau$, and since the diagram $\left( 0
\parpair{\sigma}{\tau} 1 \right)$ cannot be completed to a fork, this says
that $(X\sigma) x_0 \neq (X\tau) x_0$ for all $x_0 \in X_0$, which we
already know.  So a nondegenerate functor on $\scat{A}$ is a parallel pair
of injections with disjoint images, as claimed in~\S\ref{sec:sss}.
\end{example}

\begin{example} \label{eg:nondegen-cofork}
Let $\scat{A}$ be the category generated by objects and
arrows~\bref{eq:cofork} subject to $\rho\sigma = \rho\tau$, and consider a
functor $X: \scat{A} \go \Set$.  The nontrivial cases of \cstyle{ND1} are:
\begin{itemize}
\item $f, f' \in \{\sigma, \tau\}$: then as in
  Example~\ref{eg:nondegen-Freyd}, \cstyle{ND1} says that $X\sigma$ and
  $X\tau$ are injections with disjoint images
\item $f = f' = \rho$: the pairs of maps completing the diagram 
\[
\begin{diagdiag}
1	&		&	&		&1	\\
	&\rdTo<\rho	&	&\ldTo>\rho	&	\\
	&		&2	&		&	\\
\end{diagdiag}
\]
to a commutative square are
\[
(\id, \id), 
\ 
(\sigma, \sigma),
\ 
(\tau, \tau),
\ 
(\sigma, \tau),
\ 
(\tau, \sigma),
\]
so \cstyle{ND1} says that if $x_1, x'_1 \in X_1$ and $(X\rho) x_1 =
(X\rho) x'_1$ then $x_1 = x'_1$ (first three cases) or there exists $x_0
\in X_0$ such that $x_1 = (X\sigma) x_0$ and $x'_1 = (X\tau) x_0$ (fourth
case) or \latin{vice versa} (fifth case).  
\item $f = f' = \rho\sigma$: then \cstyle{ND1} says that $X(\rho\sigma)$ is
  injective
\item $f = \rho$, $f' = \rho\sigma$: this can be seen to be redundant.
\end{itemize}
The only nontrivial case of \cstyle{ND2} is $f = \sigma$, $f' = \tau$, and
as we saw in Example~\ref{eg:nondegen-Freyd}, this too is redundant.  So
$X$ is nondegenerate just when:
\begin{itemize}
\item $X\sigma$, $X\tau$, and $X(\rho\sigma)$ are injective
\item $X\sigma$ and $X\tau$ have disjoint images
\item if $(X\rho)x_1 = (X\rho)x'_1$ then $x_1 = x'_1$ or there exists $x_0$
  such that $\{ x_1, x'_1 \} = \{ (X\sigma)x_0, (X\tau)x_0 \}$.
\end{itemize}
An example of a nondegenerate functor on $\scat{A}$ is the diagram
\[
\{ \star \} \parpair{0}{1} [0, 1] \go S^1
\]
exhibiting the circle as an interval with its endpoints identified. 
\end{example}

\begin{example} \label{eg:nondegen-coglobular}
Let $\scat{A}$ be the category generated by objects and arrows
\[
0 \parpair{\sigma_1}{\tau_1} 1 \parpair{\sigma_2}{\tau_2} \ \cdots
\]
subject to $\sigma_{k+1} \sigma_k = \tau_{k+1} \sigma_k$ and $\sigma_{k+1}
\tau_k = \tau_{k+1} \tau_k$ for all $k\geq 1$.  A functor $\scat{A}^\op \go
\Set$ is usually called a globular set or an $\omega$-graph.  It can be
shown that a coglobular set $X: \scat{A} \go \Set$ is nondegenerate
precisely when
\begin{itemize}
\item for all $k \geq 1$, $X\sigma_k$ and $X\tau_k$ are injective
\item for all $k\geq 1$ and $x, x' \in X_k$ satisfying $(X \sigma_{k+1}) x = (X
  \tau_{k+1}) x'$, we have $x = x' \in \mr{image}(X\sigma_k) \cup
  \mr{image}(X\tau_k)$ 
\item the images of $X\sigma_1$ and $X\tau_1$ are disjoint.
\end{itemize}
For instance, the underlying coglobular set of any disk in the sense of
Joyal (\cite{Joy},~\cite{SDN}) is nondegenerate.
\end{example}

We finish this section with a diagrammatic formulation of nondegeneracy of
a module.  This will be invaluable later.

First observe that the notion of commutative diagram in a category
$\scat{A}$ can be extended to include elements of a module $M: \scat{A}
\gomod \scat{A}$.  For instance, the diagram
\[
\begin{diagram}
a_2	&\rMod^{m_2}	&a_1	&\rMod^{m_1}	&a_0	\\
\dTo<{f_2}&		&\dTo>{f_1}&		&\dTo>{f_0}\\
a'_2	&\rMod_{m'_2}	&a'_1	&\rMod_{m'_1}	&a'_0	\\
\end{diagram}
\]
is said to \demph{commute} if $m'_2 f_2 = f_1 m_2$ and $m'_1 f_1 = f_0
m_1$.  (We never attempt to compose paths containing more than one crossed
arrow $\gomod$.)  Similarly, the diagram
\[
\begin{diagram}
a_1			&\rMod^{m}	&a_0			\\
\dTo<{f_1} \dTo>{g_1}	&		&\dTo<{f_0} \dTo>{g_0}	\\
a'_1			&\rMod_{m'}	&a'_0			\\
\end{diagram}
\]
\demph{commutes serially} if $m' f_1 = f_0 m$ and $m' g_1 = g_0 m$, and the
diagram 
\[
\begin{diagram}
b	&\rMod^m	&a	&\pile{\rTo^f\\ \rTo_{f'}}	&c	\\
\end{diagram}
\]
is a \demph{fork} if $f m = f' m$.

In this language, $M: \scat{A} \gomod \scat{A}$ is nondegenerate if and
only if
\begin{description}
\item[\cstyle{ND1}]
any commutative square of solid arrows
\[
\begin{diagram}[size=2em]
	&		&b	&		&	\\
	&\ldMod(2,4)<m	&\dModget>{\!\!\!\!\!p}
				&\rdMod(2,4)>{m'}&	\\
	&		&d	&		&	\\
	&\ldGet>g	&	&\rdGet<{g'}	&	\\
a	&		&	&		&a'	\\
	&\rdTo<f	&	&\ldTo>{f'}	&	\\
	&		&c,	&		&	\\
\end{diagram}
\]
can be filled in by dotted arrows to a commutative diagram as shown, and
\item[\cstyle{ND2}] any fork $b \gobymod{m} a \parpair{f}{f'} c$
can be extended to a diagram
\[
\begin{diagram}
	&		&d			\\
	&\ruModget<p	&\dGet>e		\\
b	&\rMod_m	&a			\\
	&		&\dTo<f\dTo>{f'}	\\
	&		&c
\end{diagram}
\]
in which the triangle commutes and the right-hand column is a fork.
\end{description}

\section{Coalgebras}
\label{sec:coalgs}

We still need to prove that for any self-similarity system $(\scat{A}, M)$,
the endofunctor $M \otimes \dashbk$ of $\ftrcat{\scat{A}}{\Set}$ restricts
to an endofunctor of $\ndSet{\scat{A}}$, and similarly with $\Top$ in place
of $\Set$.  The set-theoretic case is straightforward.

\begin{propn}
\label{propn:preservationofnondegeneracy}
Let $\scat{A}$ be a small category and $M: \scat{A} \gomod \scat{A}$ a
nondegenerate module.  Then
\[
M \otimes \dashbk: \ftrcat{\scat{A}}{\Set} \go \ftrcat{\scat{A}}{\Set}
\]
preserves nondegeneracy.
\end{propn}
Nondegeneracy of $M$ is also a \emph{necessary} condition for $M \otimes
\dashbk$ to preserve nondegeneracy: consider representables.
\begin{proof}
Let $X: \scat{A} \go \Set$ be nondegenerate.  Then for any finite connected
limit $\int_i Y_i$ in $\pshf{\scat{A}}$,
\[
\int_i (Y_i \otimes M \otimes X)
\iso
\left(
\int_i (Y_i \otimes M)
\right)
\otimes X
\iso
\left(
\int_i Y_i
\right)
\otimes M \otimes X,
\]
the first isomorphism by nondegeneracy of $X$ and the second by
nondegeneracy of $M$.  So $M \otimes X$ is nondegenerate.
\done
\end{proof}

The topological case requires some preparatory lemmas.  The first concerns
$\Set$-valued functors and follows immediately from
Lemma~\ref{lemma:equalityinatensorproduct}.
\begin{lemma}[Equality in $M\otimes X$]
\label{lemma:equalityinMtensorX}
Let $\scat{A}$ be a small category, let $M: \scat{A} \gomod \scat{A}$, and
let $X \in \ndSet{\scat{A}}$.  Take module elements
\[
\begin{diagdiag}
b	&		&	&		&b'	\\
	&\rdMod<m	&	&\ldMod>{m'}	&	\\
	&		&a	&		&	\\
\end{diagdiag}
\]
and $x \in Xb$, $x' \in Xb'$.  Then $m \otimes x = m' \otimes x' \in (M
\otimes X) a$ if and only if there exist a commutative square
\[
\begin{diagdiag}
	&		&c	&		&	\\
	&\ldTo<f	&	&\rdTo>{f'}	&	\\
b	&		&	&		&b'	\\
	&\rdMod<m	&	&\ldMod>{m'}	&	\\
	&		&a	&		&	\\
\end{diagdiag}
\]
and an element $z \in Xc$ such that $fz = x$ and $f'z = x'$.
\done
\end{lemma}

\begin{lemma}[Closed quotient map]
\label{lemma:closedquotientmap}
Let $\scat{A}$ be a small category, $X: \scat{A} \go \Top$ a nondegenerate
functor, and $Y: \scat{A}^\op \go \Set$ a functor whose category of
elements is finite.  Then the quotient map
\[
q: 
\sum_a Ya \times Xa
\go
\int^a Ya \times Xa = Y \otimes X
\]
is closed. 
\end{lemma}

\begin{proof}
A subset of $Y \otimes X$ is closed just when its inverse image under $q$
is closed, so we must show that if $V$ is a closed subset of $\sum Ya
\times Xa$ then its saturation $[V] = q^{-1}qV$ is also closed.  Given $a
\in \scat{A}$ and $y \in Ya$, write $V_{a, y}$ for the intersection of $V$
with the $(a, y)$-summand $Xa$ of
\[
\sum_{(a, y) \in \elt{Y}} Xa 
\iso
\sum_{a \in \scat{A}} Ya \times Xa.
\]
Then $[V] = \bigcup_{(a, y) \in \elt{Y}} [V_{a, y}]$, so by finiteness of
$\elt{Y}$ it suffices to show that each $[V_{a, y}]$ is closed.  

Fix $(a, y) \in \elt{Y}$.  By definition,
\[
[V_{a, y}]
=
\{
(a', y', x') 
\in
\sum_{a' \in \scat{A}}
Ya' \times Xa'
\such
y' \otimes x' = y \otimes x \textrm{ for some } x \in V_{a, y}
\}.
\]
So by nondegeneracy of $X$ and Lemma~\ref{lemma:equalityinMtensorX}, $(a',
y', x') \in [V_{a, y}]$ if and only if:
\begin{condition}
there exist a span
\[
\begin{diagdiag}
	&	&b	&	&	\\
	&\ldTo<f&	&\rdTo>{f'}&	\\
a	&	&	&	&a'	\\
\end{diagdiag}
\]
in $\scat{A}$ and $z \in Xb$ such that $fz \in V_{a, y}$, $f'z = x'$, and
$yf = y'f'$,
\end{condition}
or equivalently:
\begin{condition}
there exist a span
\begin{equation}	\label{eq:spanforclosedness}
\begin{diagdiag}
	&	&(b, w)	&	&	\\
	&\ldTo<f&	&\rdTo>{f'}&	\\
(a, y)	&	&	&	&(a', y')\\
\end{diagdiag}
\end{equation}
in $\elt{Y}$ and $z \in Xb$ such that $fz \in V_{a, y}$ and $f'z = x'$.
\end{condition}
So
\[
[V_{a, y}]
=
\bigcup_{\mr{spans\ }\bref{eq:spanforclosedness}}
\{
(a', y', x') 
\such
x' \in (Xf')(Xf)^{-1} V_{a, y}
\}.
\]
But each $Xf$ is continuous and each $Xf'$ closed, so each of the sets $\{
\ldots \}$ in this union is a closed subset of the $(a', y')$-summand
$Xa'$.  Moreover, finiteness of $\elt{Y}$ guarantees that there are only
finitely many spans of the form~\bref{eq:spanforclosedness}.  So $[V_{a,
y}]$ is a finite union of closed sets, hence closed.  \done
\end{proof}

\begin{lemma}[Change of category]
\label{lemma:changeofcategory}
Let $\cat{E}$ and $\cat{E}'$ be categories with finite colimits, $F:
\cat{E} \go \cat{E}'$ a functor preserving finite colimits, $\scat{A}$ a
small category, and $M: \scat{A} \gomod \scat{A}$ a finite module.  Then
the square
\[
\begin{diagram}
\ftrcat{\scat{A}}{\cat{E}}	&\rTo^{M \otimes \dashbk}	&
\ftrcat{\scat{A}}{\cat{E}}	\\
\dTo<{F \of \dashbk}		&				&
\dTo>{F \of \dashbk}		\\
\ftrcat{\scat{A}}{\cat{E}'}	&\rTo_{M \otimes \dashbk}	&
\ftrcat{\scat{A}}{\cat{E}'},	\\
\end{diagram}
\]
commutes up to canonical isomorphism. 
\end{lemma}
\begin{proof}
Straightforward.
\done
\end{proof}

\begin{propn}
\label{propn:topologicalmoduleaction}
For all self-similarity systems $(\scat{A}, M)$, the functor
\[
M \otimes \dashbk: 
\ftrcat{\scat{A}}{\Top}
\go
\ftrcat{\scat{A}}{\Top}
\]
preserves nondegeneracy.
\end{propn}
\begin{proof}
Let $X \in \ndTop{\scat{A}}$.  The functor $U: \Top \go \Set$ preserves
colimits (being left adjoint to the indiscrete space functor), so $M
\otimes (U \of X) \iso U \of (M \otimes X)$ by
Lemma~\ref{lemma:changeofcategory}.  But $M \otimes (U \of X)$ is
nondegenerate by Proposition~\ref{propn:preservationofnondegeneracy}, so $U
\of (M \otimes X)$ is nondegenerate.

Now let $a \goby{f} a'$ be a map in $\scat{A}$, and consider the
commutative square
\[
\begin{diagram}[height=6ex]
\sum_b M(b, a) \times Xb		&
\rTo^{\sum f_* \times 1}		&
\sum_b M(b, a') \times Xb		\\
\dQt<{q_a}				&
					&
\dQt>{q_{a'}}				\\
(M \otimes X) a				&
\rTo_{(M \otimes X) f}			&
(M \otimes X) a'.			\\
\end{diagram}
\]
The map $\sum f_* \times 1$ is closed because each set $M(b, a)$ is finite.
The map $q_{a'}$ is closed by Lemma~\ref{lemma:closedquotientmap} and
finiteness of $M$.  So $((M \otimes X) f) \of q_a$ is closed; but $q_a$ is
a continuous surjection, so $(M \otimes X) f$ is closed.  \done
\end{proof}

Any $M$-coalgebra $(X, \xi)$ in $\Top$ has an underlying $M$-coalgebra in
$\Set$, since $U \of X: \scat{A} \go \Set$ is nondegenerate and there is a
natural transformation
\[
U \xi: 
U \of X 
\go
U \of (M \otimes X)
\iso
M \otimes (U \of X).
\]
This defines a functor
\[
U_*: 
\Coalg{M}{\Top}
\go
\Coalg{M}{\Set}
\]
where the domain and codomain are the categories of $M$-coalgebras in
$\Top$ and $\Set$ respectively.

\begin{propn}[$\Top$ \vs\ $\Set$]
\label{propn:coalgebrasinTopandSet}
Let $(\scat{A}, M)$ be a self-similarity system.  The forgetful functor
\[
U_*: 
\Coalg{M}{\Top}
\go
\Coalg{M}{\Set}
\]
has a left adjoint, and if $(I, \iota)$ is a universal solution in $\Top$
then $U_*(I, \iota)$ is a universal solution in $\Set$. 
\end{propn}
Conversely, we will see that any universal solution in $\Set$ can be
equipped with a topology that makes it the universal solution in $\Top$.

\begin{proof}
Let $D$ be the left adjoint to $U: \Top \go \Set$, assigning to each set
the corresponding discrete space.  Then $D$ preserves colimits, so commutes
with $M \otimes \dashbk$ (Lemma~\ref{lemma:changeofcategory}).  Moreover,
if $X: \scat{A} \go \Set$ is nondegenerate then so is $D \of X: \scat{A}
\go \Top$.  Hence $D$ induces a functor
\[
D_*: 
\Coalg{M}{\Set} 
\go 
\Coalg{M}{\Top}.
\]
For purely formal reasons, the adjunction $D \ladj U$ induces an adjunction
$D_* \ladj U_*$.  The statement on universal solutions follows from the
fact that right adjoints preserve terminal objects.
\done 
\end{proof}

%% file: univsoln.tex
\section{The universal solution}
\label{sec:univ-soln}

In this section I construct the universal solution in $\Set$ and in $\Top$
of any given self-similarity system, assuming that the system satisfies a
certain solvability condition \So.  In the Appendix I show that this
sufficient condition \So\ is also necessary.  The construction therefore
gives the universal solution whenever one exists.

Condition \So\ on a self-similarity system $(\scat{A}, M)$ is:
\begin{description}
\item[\cstyle{S1}]
given any commutative diagram
\[
\begin{diagram}
\cdots		&\rMod^{m_3}	&a_2		&\rMod^{m_2}	&
a_1		&\rMod^{m_1}	&a_0		\\
		&		&\dTo>{f_2}	&		&
\dTo>{f_1}	&		&\dTo>{f_0}	\\
\cdots		&\rMod^{p_3}	&b_2		&\rMod^{p_2}	&
b_1		&\rMod^{p_1}	&b_0		\\
		&		&\uTo>{f'_2}	&		&
\uTo>{f'_1}	&		&\uTo>{f'_0}	\\
\cdots		&\rMod_{m'_3}	&a'_2		&\rMod_{m'_2}	&
a'_1		&\rMod_{m'_1}	&a'_0,		\\
\end{diagram}
\]
there exists a commutative square
\[
\begin{diagdiag}
	&	&a_0	&		&	\\
	&\ruTo	&	&\rdTo>{f_0}	&	\\
\cdot	&	&	&		&b_0	\\
	&\rdTo	&	&\ruTo>{f'_0}	&	\\
	&	&a'_0	&		&	\\
\end{diagdiag}
\]
in $\scat{A}$, and
\item[\cstyle{S2}] given any serially commutative diagram
\[
\begin{diagram}
\cdots		&\rMod^{m_3}	&a_2			&\rMod^{m_2}	&
a_1			&\rMod^{m_1}	&a_0			\\
		&		&\dTo<{f_2} \dTo>{f'_2}	&		&
\dTo<{f_1} \dTo>{f'_1} 	&		&\dTo<{f_0} \dTo>{f'_0}	\\
\cdots		&\rMod_{p_3}	&b_2			&\rMod_{p_2}	&
b_1			&\rMod_{p_1}	&b_0,		\\
\end{diagram}
\]
there exists a fork $\cdot \go a_0 \parpair{f_0}{f'_0} b_0$ in
$\scat{A}$. 
\end{description}

\begin{example}
For any small category $\scat{A}$ there is a module $M: \scat{A} \gomod
\scat{A}$ defined by $M(b, a) = \scat{A}(b, a)$, and $(\scat{A}, M)$ is a
self-similarity system as long as $\sum_b \scat{A}(b, a)$ is finite for
each $a \in \scat{A}$.  Condition \So\ says that $\scat{A}$ is
componentwise cofiltered; so, for instance, the self-similarity system
obtained by taking $\scat{A} = ( 0 \parpairu 1 )$ has no universal
solution.  If $\scat{A}$ \emph{is} componentwise cofiltered then the
universal solution is the functor $\scat{A} \go \Top$ constant at the
one-point space, with its unique coalgebra structure.
\end{example}

Here is the construction of the universal solution.  The proofs that it
works are in~\S\ref{sec:Set-proofs} (for $\Set$) and~\S\ref{sec:Top-proofs}
(for $\Top$).

Let $(\scat{A}, M)$ be a self-similarity system.  For each $a \in
\scat{A}$, there is a category $\catI a$ in which an object is an infinite
sequence
\begin{equation}	\label{eq:typical-sequence}
\cdots \gobymod{m_3} a_2 \gobymod{m_2} a_1 \gobymod{m_1} a_0 = a
\end{equation}
and a map $(a_\blob, m_\blob) \go (a'_\blob, m'_\blob)$ is a commutative
diagram
\begin{equation}	\label{eq:typical-ladder}
\begin{diagram}
\cdots		&\rMod^{m_3}	&a_2		&\rMod^{m_2}	&
a_1		&\rMod^{m_1}	&a_0 = a		\\
		&		&\dTo>{f_2}	&		&
\dTo>{f_1}	&		&\dTo>{f_0 = 1_a}	\\
\cdots		&\rMod_{m'_3}	&a'_2		&\rMod_{m'_2}	&
a'_1		&\rMod_{m'_1}	&a'_0 = a.		\\
\end{diagram}
\end{equation}
Moreover, each map $f: a \go a'$ in $\scat{A}$ induces a functor $\catI f:
\catI a \go \catI a'$, sending an object~\bref{eq:typical-sequence} of
$\catI a$ to the object
\[
\cdots \gobymod{m_3} a_2 \gobymod{m_2} a_1 \gobymod{f m_1} a'
\]
of $\catI a'$.  This defines a functor $\catI : \scat{A} \go \Cat$. 

Write $\Pi_0: \Cat \go \Set$ for the functor sending a small category to
its set of connected-components, and put $I = \Pi_0 \catI : \scat{A} \go
\Set$; thus, $Ia$ is a set of equivalence classes of
diagrams~\bref{eq:typical-sequence}.  In~\S\ref{sec:Set-proofs} we will see
that if $(\scat{A}, M)$ satisfies~\So\ then $I$ is nondegenerate.

\begin{warning} \label{warning:fin}
$Ia$ is \emph{not} the limit of finite approximations.  Precisely, let
$\catI_n a$ be the category whose objects are diagrams of the form
\begin{equation} \label{eq:finitesequence}
a_n \gobymod{m_n} \cdots \gobymod{m_1} a_0 = a
\end{equation}
and whose arrows are commutative diagrams, and let $I_n a$ be the set of
connected-components of $\catI_n a$: then $\catI a$ is the limit of the
$(\catI_n a)$s, but $Ia$ is not in general the limit of the $(I_n a)$s.
See the last paragraph of this section for an example.
\end{warning}

There is an $M$-coalgebra structure $\iota$ on $I$ defined by
\[
\iota_a
\bigcompt{
\cdots \gobymod{m_3} a_2 \gobymod{m_2} a_1 \gobymod{m_1} a
}
=
(
a_1 \gobymod{m_1} a
)
\otimes
\bigcompt{
\cdots \gobymod{m_3} a_2 \gobymod{m_2} a_1
}
\]
where $\compt{\emptybk}$ denotes connected-component.  To see that this is
a valid definition, note that given a map~\bref{eq:typical-ladder} in
$\catI a$, there is an equality $m_1 = m'_1 f_1$ and a map
\[
\begin{diagram}
\cdots	&\rMod^{m_3}	&a_2		&\rMod^{f_1 m_2}	&a'_1	\\
	&		&\dTo>{f_2}	&			&\dTo>1	\\
\cdots	&\rMod_{m'_3}	&a'_2		&\rMod_{m'_2}		&a'_1	\\
\end{diagram}
\]
in $\catI a'_1$, so 
\begin{eqnarray*}
m_1
\otimes
\bigcompt{ 
\cdots \gobymod{m_3} a_2 \gobymod{m_2} a_1
}	&
=	&
m'_1
\otimes
\bigcompt{ 
\cdots \gobymod{m_3} a_2 \gobymod{f_1 m_2} a'_1
}	\\
	&
=	&
m'_1
\otimes
\bigcompt{ 
\cdots \gobymod{m'_3} a'_2 \gobymod{m'_2} a'_1
}.
\end{eqnarray*}
Naturality of $\iota$ is easily checked.  So if \So\ holds then $(I,
\iota)$ is an $M$-coalgebra, and it is in fact the terminal one, that is,
the universal solution in $\Set$
(Theorem~\ref{thm:universalsolutioninSet}).

Now we construct the topology.  For each $a \in \scat{A}$, $n \in \nat$,
and diagram of the form~\bref{eq:finitesequence}, there is a subset of
$Ia$ consisting of all those $t$ such that
\[
t = 
\bigcompt{ 
\cdots 
\gobymod{m_{n+2}} a_{n+1} \gobymod{m_{n+1}} a_n \gobymod{m_n}
\cdots
\gobymod{m_1} a_0 = a
}
\]
for some $m_{n+1}, a_{n+1}, m_{n+2}, \ldots$.  Generate a topology on $Ia$
by taking each such subset to be closed.  It is not obvious that the maps
$If: Ia \go Ia'$ are continuous or closed, or that the maps $\iota_a: Ia
\go (M\otimes I)a$ are continuous.  Nevertheless, if~\So\ holds then
they are, and so $(I, \iota)$ is an $M$-coalgebra in $\Top$.
Theorem~\ref{thm:universalsolutioninTop} says that it is in fact the
universal solution in $\Top$.

(It is not hard to see that each of these basic closed subsets must
be closed in \emph{any} topology on 
$(I, \iota)$~\cite[\lemmafixedpointcomponents]{SS2}.  So this is the 
coarsest possible topology.)

Let us see how all of this works in the case of the Freyd self-similarity
system $(\scat{A}, M)$.

First, condition~\So\ holds.  For \cstyle{S1}, the only diagram
\[
\begin{diagdiag}
a_0	&		&	\\
	&\rdTo>{f_0}	&	\\
	&		&b_0	\\
	&\ruTo>{f'_0}	&	\\
a'_0	&		&	\\
\end{diagdiag}
\]
in $\scat{A}$ that cannot be completed to a commutative square is (up to
symmetry) that in which $f_0 = \sigma$ and $f'_0 = \tau$, and then there is
no infinite commutative diagram as in \cstyle{S1}: indeed, there is not
even a commutative diagram of the form
\[
\begin{diagram}
\cdot	&\rMod	&0		\\
\dTo	&	&\dTo>\sigma	\\
\cdot	&\rMod	&1		\\
\uTo	&	&\uTo>\tau	\\
\cdot	&\rMod	&0.		\\
\end{diagram}
\]
Similarly, for \cstyle{S2}, the only parallel pair of arrows in $\scat{A}$
that cannot be completed to a fork is $0 \parpair{\sigma}{\tau} 1$, and
there is no serially commutative diagram of the form
\[
\begin{diagram}
\cdot		&\rMod	&0			\\
\dTo\dTo	&	&\dTo<\sigma \dTo>\tau	\\
\cdot		&\rMod	&1.			\\
\end{diagram}
\]

The universal solution $(I, \iota)$ has $I1 = [0, 1]$, so according to the
construction, an element of $[0, 1]$ is an equivalence class of diagrams
\[
\cdots \gobymod{m_3} a_2 \gobymod{m_2} a_1 \gobymod{m_1} 1.
\]
If each $a_n$ is $1$ then each $m_n$ is either $[0, \half]$ or $[\half, 1]$
and the diagram is just a binary expansion: for instance,
\[
\cdots \gobymod{[\half, 1]} 1 \gobymod{[0, \half]} 1 
\gobymod{[\half, 1]} 1 \gobymod{[0, \half]} 1
\gobymod{[\half, 1]} 1
\]
corresponds to $0.10101\ldots$, representing $\frc{2}{3} \in [0, 1]$.
Otherwise, the diagram is of the form
\[
\cdots 
\gobymod{\id} 0 \gobymod{\id} 0
\gobymod{m_{n+1}} 1 \gobymod{m_n}
\cdots
\gobymod{m_1} 1
\]
where $m_1, \ldots, m_n \in \{ [0, \half], [\half, 1] \}$ and $m_{n+1} \in
\{ 0, \half, 1 \}$: for instance,
\[
\cdots 
\gobymod{\id} 0 \gobymod{\id} 0
\gobymod{\half} 1 \gobymod{[\half, 1]} 1
\gobymod{[\half, 1]} 1 \gobymod{[0, \half]} 1.
\]
To see which element $t$ of $[0, 1]$ this represents, reason as follows:
the $[0, \half]$ says that $t \in [0, \half]$; the right-hand copy of
$[\half, 1]$ says that $t$ is in the upper half of $[0, \half]$, that is,
in $[\frc{1}{4}, \half]$; the left-hand copy of $[\half, 1]$ says that $t$
is in the upper half of $[\frc{1}{4}, \half]$, that is, in $[\frc{3}{8},
\half]$; then the $\half$ says that $t$ is the midpoint of $[\frc{3}{8},
\half]$, that is, $t = \frc{7}{16}$.

An element of $[0, 1]$ has at most two binary expansions, but may have
infinitely many representations in $\catI 1$.  For instance, the
representations of $\half$ are
\begin{eqnarray}
\cdots \gobymod{[\half, 1]} 1 \gobymod{[\half, 1]} 1 
\gobymod{[0, \half]} 1,
\label{eq:half-1}	\\
\cdots \gobymod{[0, \half]} 1 \gobymod{[0, \half]} 1 
\gobymod{[\half, 1]} 1, 
\label{eq:half-2}	\\
\cdots \gobymod{\id} 0 \gobymod{\id} 0 \gobymod{\half} 1,
\label{eq:half-3}  
\end{eqnarray}
and for any $n\in\nat$,
\begin{eqnarray}
\cdots \gobymod{\id} 0 \gobymod{1} 1
\gobymod{[\half, 1]} \cdots \gobymod{[\half, 1]} 1
\gobymod{[0, \half]} 1,
\label{eq:half-4}	\\
\cdots \gobymod{\id} 0 \gobymod{0} 1
\gobymod{[0, \half]} \cdots \gobymod{[0, \half]} 1
\gobymod{[\half, 1]} 1\ 
\label{eq:half-5}	
\end{eqnarray}
with $n$ copies of $[\half, 1]$ and $[0, \half]$ respectively.

The construction says that two objects of $\catI 1$ represent the same
element of $[0, 1]$ if and only if they are in the same
connected-component.  So, for instance, each of
\bref{eq:half-1}--\bref{eq:half-5} should be in the same component; the
connected diagram
\[
\begin{diagram}
\cdots	&\rMod^{[\half, 1]}	&1		&
\rMod^{[\half, 1]}	&1		&
\rMod^{[0, \half]}	&1	\\
	&			&\uTo>\tau	&
			&\uTo>\tau	&
			&\uTo>1	\\
\cdots	&\rMod^{\id}		&0		&
\rMod^{\id}		&0		&
\rMod^{\half}		&1	\\
	&			&\dTo>\sigma	&
			&\dTo>\sigma	&
			&\dTo>1	\\
\cdots	&\rMod_{[0, \half]}	&1		&
\rMod_{[0, \half]}	&1		&
\rMod_{[\half, 1]}	&1	\\
\end{diagram}
\]
shows that \bref{eq:half-1}--\bref{eq:half-3} are, and the others are left
to the reader.  Observe%
\label{p:obs-finitary-move}
in general that any two objects of $\catI a$ of the form
\begin{eqnarray}
\cdots 
\gobymod{m'_{n+2}} a'_{n+1} \gobymod{m'_{n+1}} a'_n 
\gobymod{m_n f} a_{n-1} \gobymod{m_{n-1}} 
\cdots
\gobymod{m_1} a_0 = a,
\label{eq:finitary-right}
\\
\cdots 
\gobymod{m'_{n+2}} a'_{n+1} \gobymod{fm'_{n+1}} a_n 
\gobymod{m_n} a_{n-1} \gobymod{m_{n-1}} 
\cdots
\gobymod{m_1} a_0 = a\ 
\label{eq:finitary-left}
\end{eqnarray}
(where $f: a'_n \go a_n$) are in the same connected-component, because
there is a map
\[
( \ldots, 1_{a'_{n+1}}, f, 1_{a_{n-1}}, \ldots, 1_{a_0} )
\]
from the top row to the bottom.

The constructed topology on $[0, 1]$ is generated by taking as closed
all subsets of the form $[k/2^n, l/2^n]$ where $k, l, n \in \nat$ and $l
\in \{k, k + 1\}$.  This is exactly the Euclidean topology.

Finally, this example shows that $Ia$ need not be the limit of finite
approximations (Warning~\ref{warning:fin}).  It is not hard to show that
for each $n$, the category $\catI_n 1$ is connected: so each $I_n 1$ is a
one-element set, and $I1 = [0, 1]$ is plainly not the limit of the $(I_n
1)$s.

\section{Set-theoretic proofs}
\label{sec:Set-proofs}

Fix a self-similarity system $(\scat{A}, M)$ satisfying \So.  In this
section I prove that the functor $I: \scat{A} \go \Set$ is nondegenerate
and that $(I, \iota)$ is the universal solution.

Nondegeneracy of $I$ will follow from a kind of nondegeneracy property of
$\catI: \scat{A} \go \Cat$.  Any $\Cat$-valued functor $\cat{X}: \scat{B}
\go \Cat$ has a \demph{category of elements} $\elt{\cat{X}}$, in which an
object is a pair $(b, x)$ with $b \in \scat{B}$ and $x \in \cat{X}b$ and an
arrow $(b, x) \go (b', x')$ is a pair $(g, \xi)$ with $g: b \go b'$ in
$\scat{B}$ and $\xi: (\cat{X}g) x \go x'$ in $\cat{X}b'$.  This is related
to the notion of the category of elements of a $\Set$-valued functor $X:
\scat{B} \go \Set$ (page~\pageref{p:cat-elts}) by the isomorphism $\elt{X}
\iso \elt{D \of X}$, where $D: \Set \go \Cat$ is the functor assigning to
each set the corresponding discrete category.  The nondegeneracy property
mentioned is that $\elt{\catI}$ is componentwise cofiltered (that is,
$\elt{\catI}^\op$ is componentwise filtered, in the sense defined after
Corollary~\ref{cor:componentwisefilteredcategories}).

This is proved by a finiteness argument.  Notation: if $\cat{L}$ is the
limit of a diagram
\[
\cdots \go \cat{L}_3 \go \cat{L}_2 \go \cat{L}_1
\]
in some category, I write $\pr_n$ for both the projection $\cat{L} \go
\cat{L}_n$ and the given map $\cat{L}_m \go \cat{L}_n$ for any $m \geq n$.

\begin{lemma}[K\"onig~\cite{Koen}]
\label{lemma:Koenig}
The limit in $\Set$ of a diagram
\[
\cdots \go \cat{F}_3 \go \cat{F}_2 \go \cat{F}_1
\]
of finite nonempty sets is nonempty.  More precisely, for any sequence
$(F_n)_{n \geq 1}$ with $F_n \in \cat{F}_n$ there exists an element $G$ of
the limit such that
\[
\forall r \geq 1, 
\ 
\exists n \geq r:
\ 
\pr_r (F_n) = \pr_r (G).
\]
\end{lemma}

\paragraph*{Remark}  The first sentence is a special case of the fact
that a componentwise cofiltered limit of nonempty compact Hausdorff spaces
is nonempty (compare~\cite[I.9.6]{Bou}).  
 
\begin{proof}
Take a sequence $(F_n)_{n \geq 1}$ with $F_n \in \cat{F}_n$.  We define,
for each $r \geq 1$, an infinite subset $\nat_r$ of $\nat$ and an element
$G_r \in \cat{F}_r$ such that
\begin{itemize}
\item for all $r \geq 1$, $\nat_r \sub \nat_{r-1} \cap \{ r, r+1, \ldots
\}$ (writing $\nat_0 = \nat$)
\item for all $r \geq 1$ and $n \in \nat_r$, $\pr_r (F_n) = G_r$.
\end{itemize}
Suppose inductively that $r \geq 1$ and $\nat_{r-1}$ is defined.  As $n$
runs through the infinite set $\nat_{r-1} \cap \{ r, r+1, \ldots \}$,
$\pr_r (F_n)$ takes values in the finite set $\cat{F}_r$, so takes some
value $G_r \in \cat{F}_r$ infinitely often.  Putting
\[
\nat_r 
= 
\{ 
n \in \nat_{r-1} \cap \{ r, r+1, \ldots \}
\such
\pr_r (F_n) = G_r 
\}
\]
completes the induction.

For each $r \geq 1$ we have $G_r = \pr_r (G_{r+1})$, since we may
choose $n \in \nat_{r+1}$ and then
\[
\pr_r (G_{r+1})
=
\pr_r (\pr_{r+1} (F_n))
=
\pr_r (F_n)
=
G_r.
\]
So there is a unique element $G$ of the limit such that $\pr_r (G) = G_r$
for all $r \geq 1$.  Given $r \geq 1$, we may choose $n \in \nat_r$, and
then $n \geq r$ and $\pr_r (F_n) = G_r = \pr_r (G)$ as required. 
\done
\end{proof}

This will be applied as follows.  Suppose we have a limit $\cat{L}$ of
categories
\[
\cdots \go \cat{L}_3 \go \cat{L}_2 \go \cat{L}_1,
\]
categories $\scat{J} \sub \scat{K}$, and a diagram $D$ of shape $\scat{J}$
in $\cat{L}$, and suppose we are interested in extending $D$ to a diagram
of shape $\scat{K}$.  If we can do so then we can certainly extend $\pr_n
\of D: \scat{J} \go \cat{L}_n$ to $\scat{K}$ for all $n$; the converse does
not hold, because we cannot necessarily choose the extensions $\scat{K} \go
\cat{L}_n$ in a coherent way.  But the following lemma says that we can do
it if for each $n$ there are only finitely many choices of extension.

\begin{lemma}[Factorization]
\label{lemma:factorization}
Let $\cat{L}$ be the limit of a diagram
\[
\cdots \go \cat{L}_3 \go \cat{L}_2 \go \cat{L}_1
\]
in some category $\cat{E}$.  Let $P: \scat{J} \go \scat{K}$ and $D:
\scat{J} \go \cat{L}$ be maps in $\cat{E}$ such that for each $n \geq 1$,
the set of factorizations of $\pr_n \of D$ through $P$ is finite and
nonempty.  Then $D$ factors through $P$.
\end{lemma}

\begin{proof}
Apply Lemma~\ref{lemma:Koenig} with $\cat{F}_n = \{ F \in \cat{E}(\scat{K},
\cat{L}_n) \such F \of P = \pr_n \of D \}$.
\done 
\end{proof}

For each $n\in\nat$ we have a functor $\catI_n: \scat{A} \go \Cat$
(see~\ref{warning:fin}).  The evident projections make $\elt{\catI}$ the
limit in $\Cat$ of
\[
\cdots \go \elt{\catI_2} \go \elt{\catI_1}.
\]

\begin{propn}
\label{propn:catInondegenerate}
$\elt{\catI}$ is componentwise cofiltered.
\end{propn}
\begin{proof}
We have to prove that every diagram $\cdot \go \cdot \og \cdot$ in
$\elt{\catI}$ can be completed to a commutative square and that every
parallel pair $\cdot \parpairu \cdot$ can be completed to a fork.  The two
cases are very similar, so I just do the first.

Take a diagram
\begin{equation} \label{eq:corner-in-I}
\begin{diagram}
\cdots		&\rMod^{m_3}	&a_2		&\rMod^{m_2}	&
a_1		&\rMod^{m_1}	&a_0		\\
		&		&\dTo>{f_2}	&		&
\dTo>{f_1}	&		&\dTo>{f_0}	\\
\cdots		&\rMod^{p_3}	&b_2		&\rMod^{p_2}	&
b_1		&\rMod^{p_1}	&b_0		\\
		&		&\uTo>{f'_2}	&		&
\uTo>{f'_1}	&		&\uTo>{f'_0}	\\
\cdots		&\rMod_{m'_3}	&a'_2		&\rMod_{m'_2}	&
a'_1		&\rMod_{m'_1}	&a'_0		\\
\end{diagram}
\end{equation}
of shape $\littlepullback$ in $\elt{\catI}$.  We apply
Lemma~\ref{lemma:factorization} where $P$ is the inclusion
\[
\scat{J}
=
\left(
\begin{diagdiag}
	&	&\	&	&	\\
	&	&\	&	&	\\
\cdot	&	&	&	&\cdot	\\
	&\rdTo	&	&\ldTo	&	\\
	&	&\cdot	&	&	\\
\end{diagdiag}
\right)
\rIncl
\left(
\begin{diagdiag}
	&	&\cdot	&	&	\\
	&\ldTo	&	&\rdTo	&	\\
\cdot	&	&	&	&\cdot	\\
	&\rdTo	&	&\ldTo	&	\\
	&	&\cdot	&	&	\\
\end{diagdiag}
\right)
=
\scat{K}
\]
(the square in $\scat{K}$ being commutative), $\cat{L}_n =
\elt{\catI_n}$, $\cat{L} = \elt{\catI}$, and $D$ is the
diagram~\bref{eq:corner-in-I} with its rightmost block removed.  The
hypothesis of Lemma~\ref{lemma:factorization} is that for each $n\geq 1$,
the set of diagrams of the form
\begin{equation} \label{eq:factorization-chunk}
\begin{diagram}
a_n		&\rMod^{m_n}	&a_{n-1}	&\rMod^{m_{n-1}}&
\	&\cdots	&\		&\rMod^{m_2}	&a_1		\\
\uTo<{g_n}	&		&\uTo<{g_{n-1}}	&		&
	&	&		&		&\uTo>{g_1}	\\
c_n		&\rMod^{q_n}	&c_{n-1}	&\rMod^{q_{n-1}}&
\	&\cdots	&\ 		&\rMod^{q_2}	&c_1		\\
\dTo<{g'_n}	&		&\dTo<{g'_{n-1}}&		&
	&	&		&		&\dTo>{g'_1}	\\
a'_n		&\rMod_{m'_n}	&a'_{n-1}	&\rMod_{m'_{n-1}}&
\	&\cdots	&\ 		&\rMod_{m'_2}	&a'_1		\\
\end{diagram}
\end{equation}
satisfying $f_1 g_1 = f'_1 g'_1, \ldots, f_n g_n = f'_n g'_n$ is nonempty
and finite.   

Finiteness follows from finiteness of $M$ (and the absence of the rightmost
block).  For nonemptiness, let $n\geq 1$.  Then \cstyle{S1} implies that
there exist $c_n$, $g_n$, and $g'_n$ making
\[
\begin{diagram}
	&		&a_n	&		&\rMod^{m_n}	&
	&a_{n-1}	&			&	\\
	&\ruTo<{g_n}	&	&\rdTo>{f_n}	&		&
\	&		&\rdTo>{f_{n-1}}	&	\\
c_n	&		&	&		&b_n		&
	&\rMod^{p_n}	&			&b_{n-1}\\
	&\rdTo<{g'_n}	&	&\ruTo>{f'_n}	&		&
	&		&\ruTo>{f'_{n-1}}	&	\\
	&		&a'_n	&		&\rMod_{m'_n}	&
	&a'_{n-1}	&			&	\\
\end{diagram}
\]
commute, and then nondegeneracy of $M$ (condition~\cstyle{ND1} at the end
of~\S\ref{sec:nondegen}) implies that the outside of this diagram can also
be filled in as
\[
\begin{diagram}
	&		&a_n	&		&\rMod^{m_n}	&
	&a_{n-1}	&			&	\\
	&\ruTo<{g_n}	&	&\ 		&		&
\ruTo<{g_{n-1}}&	&\rdTo>{f_{n-1}}	&	\\
c_n	&		&\rMod^{q_n}&		&c_{n-1}	&
	&		&			&b_{n-1}\\
	&\rdTo<{g'_n}	&	&		&		&
\rdTo<{g'_{n-1}}&	&\ruTo>{f'_{n-1}}	&	\\
	&		&a'_n	&		&\rMod_{m'_n}	&
	&a'_{n-1}.	&			&	\\
\end{diagram}
\]
Repeating this argument $(n-2)$ times gives a
diagram~\bref{eq:factorization-chunk}, as required.

So the hypothesis of Lemma~\ref{lemma:factorization} holds, and $D$ can be
completed to a commutative square in $\elt{\catI}$.  Using the
diagram-filling argument one more time shows that~\bref{eq:corner-in-I} can
be, too.  \done
\end{proof}

The next few results show that for general reasons, $\elt{\catI}$ being
componentwise cofiltered implies that each $\catI a$ is too and that $I:
\scat{A} \go \Set$ is nondegenerate.

\begin{lemma}
\label{lemma:individualcategoriescomponentwisecofiltered}
Let $\cat{J}: \scat{B} \go \Cat$ be a functor on a small category
$\scat{B}$.  If $\elt{\cat{J}}$ is componentwise cofiltered then
$\cat{J}a$ is componentwise cofiltered for each $a \in \scat{B}$.
\end{lemma}

\begin{proof}
We have to prove that every diagram $\cdot \go \cdot \og \cdot$ in
$\cat{J}a$ can be completed to a commutative square and that every parallel
pair $\cdot \parpairu \cdot$ can be completed to a fork.  Again I just do
the first case; the second is similar.

Take a diagram
\[
\begin{diagdiag}
\omega	&		&	&		&\omega'	\\
	&\rdTo<\phi	&	&\ldTo>{\phi'}	&		\\
	&		&\chi	&		&		\\
\end{diagdiag}
\]
in $\cat{J}a$.  Then there is a commutative square
\[
\begin{diagdiag}
	&		&(b, \zeta)&		&		\\
	&\ldTo<{(g, \gamma)}&	&\rdTo>{(g', \gamma')}&		\\
(a, \omega)&		&	&		&(a, \omega')	\\
	&\rdTo<{(1, \phi)}&	&\ldTo>{(1, \phi')}&		\\
	&		&(a, \chi)&		&		\\
\end{diagdiag}
\]
in $\elt{\cat{J}}$.  Commutativity says that $g = g'$ and that the square 
\[
\begin{diagdiag}
	&		&(\cat{J}g) \zeta	&		&	\\
	&\ldTo<\gamma	&			&\rdTo>{\gamma'}&	\\
\omega	&		&			&		&\omega'\\
	&\rdTo<\phi	&			&\ldTo>{\phi'}	&	\\
	&		&\chi			&		&	\\
\end{diagdiag}
\]
in $\cat{J}a$ commutes, as required.
\done
\end{proof}

\begin{propn}
\label{propn:Iacomponentwisecofiltered}
$\catI a$ is componentwise cofiltered for each $a \in \scat{A}$.  
\done
\end{propn}

\begin{lemma}
\label{lemma:connectedcomponents}
Let $\cat{J}: \scat{B} \go \Cat$ be a functor on a small category
$\scat{B}$.  If $\elt{\cat{J}}$ is componentwise cofiltered then so is
$\elt{\Pi_0 \cat{J}}$.
\end{lemma}
\begin{proof}
Once again the proof splits into two similar cases.  For variety, I do the
second: that every diagram
$
(a, \compt{\omega}) \parpair{f}{f'} (b, \compt{\chi})
$
in $\elt{\Pi_0 \cat{J}}$ extends to a fork.  

Since $\compt{(\cat{J}f) \omega} = \compt{\chi} = \compt{(\cat{J}f')
\omega}$, Lemmas~\ref{lemma:connectednessbyspans}
and~\ref{lemma:individualcategoriescomponentwisecofiltered} imply that
there exists a span
\[
\begin{diagdiag}
		&	&\xi	&	&		\\
		&\ldTo<\delta&	&\rdTo>{\delta'}&	\\
(\cat{J}f)\omega&	&	&	&(\cat{J}f')\omega\\
\end{diagdiag}
\]
in $\cat{J}b$.  We therefore have a finite connected diagram (solid arrows)
\[
\begin{diagdiag}
			&			&
(c, \zeta)		&
			&			\\
			&\ldGet<{(g, \gamma)}	&
			&
\rdGet>{(fg = f'g, \theta)}&			\\
(a, \omega)		&			&
			&
			&(b, \xi)		\\
			&\rdTo(4,4)>{\rdlabel{(f', 1)}}	&
			&
\ldTo(4,4)<{\ldlabel{(1, \delta)}}&		\\
\dTo<{(f, 1)}		&			&
			&
			&\dTo>{(1, \delta')}	\\
			&			&
			&
			&			\\
(b, (\cat{J}f)\omega)	&			&
			&
			&(b, (\cat{J}f')\omega)\\
\end{diagdiag}
\]
in $\elt{\cat{J}}$, so by hypothesis there exists a dotted commutative
diagram, giving a fork
\[
(c, \compt{\zeta})
\goby{g}
(a, \compt{\omega})
\parpair{f}{f'}
(b, \compt{\chi})
\]
in $\elt{\Pi_0 \cat{J}}$.  
\done
\end{proof}

\begin{propn}[$I$ nondegenerate]
\label{propn:setInondegenerate}
$I: \scat{A} \go \Set$ is nondegenerate.  
\done
\end{propn}

Hence $(I, \iota)$ is an $M$-coalgebra.  By Lambek's Lemma, a necessary
condition for it to be the universal solution is that $\iota$ is an
isomorphism, and we can prove this immediately.

\begin{propn}[$I$ is a fixed point]
\label{propn:Iisafixedpoint}
$\iota: I \go M \otimes I$ is an isomorphism.
\end{propn}

\begin{proof}
It is enough to show that $\iota_a: Ia \go (M \otimes I)a$ is bijective for
each $a \in \scat{A}$.  Certainly $\iota_a$ is surjective.  For
injectivity, suppose that
\[
\iota_a
\bigcompt{
\cdots \gobymod{m_2} a_1 \gobymod{m_1} a
}
=
\iota_a
\bigcompt{
\cdots \gobymod{m'_2} a'_1 \gobymod{m'_1} a
},
\]
that is,
\[
m_1
\otimes
\bigcompt{
\cdots \gobymod{m_2} a_1
}
=
m'_1
\otimes
\bigcompt{
\cdots \gobymod{m'_2} a'_1
}.
\]
By Lemma~\ref{lemma:equalityinMtensorX} and nondegeneracy of $I$,
there exist a commutative square
\[
\begin{diagdiag}
	&		&b	&		&	\\
	&\ldTo<f	&	&\rdTo>{f'}	&	\\
a_1	&		&	&		&a'_1	\\
	&\rdMod<{m_1}	&	&\ldMod>{m'_1}	&	\\
	&		&a	&		&	\\
\end{diagdiag}
\]
and an element $\compt{\cdots \gobymod{p_2} b_1 \gobymod{p_1} b} \in Ib$
such that
\begin{eqnarray*}
\bigcompt{
\cdots \gobymod{p_2} b_1 \gobymod{f p_1} a_1
}	&
=	&
\bigcompt{
\cdots \gobymod{m_3} a_2 \gobymod{m_2} a_1
},	\\
\bigcompt{
\cdots \gobymod{p_2} b_1 \gobymod{f' p_1} a'_1
}	&
=	&
\bigcompt{
\cdots \gobymod{m'_3} a'_2 \gobymod{m'_2} a'_1
}.
\nonumber	
\end{eqnarray*}
Then 
\begin{eqnarray*}
\bigcompt{
\cdots \gobymod{m_3} a_2 \gobymod{m_2} a_1 \gobymod{m_1} a
}	&
=	&
\bigcompt{
\cdots \gobymod{p_2} b_1 \gobymod{f p_1} a_1 \gobymod{m_1} a
}	\\
	&
=	&
\bigcompt{
\cdots \gobymod{p_2} b_1 \gobymod{p_1} b \gobymod{m_1 f} a
},
\end{eqnarray*}
using the observation at~\bref{eq:finitary-right}
and~\bref{eq:finitary-left} (page~\pageref{eq:finitary-right}).  But $m_1 f
= m'_1 f'$, so by symmetry of argument,
\[
\bigcompt{
\cdots \gobymod{m_2} a_1 \gobymod{m_1} a
}
=
\bigcompt{
\cdots \gobymod{m'_2} a'_1 \gobymod{m'_1} a
}
\]
as required.
\done
\end{proof}

Here is the key concept for the rest of the proof.  Let $(X, \xi)$ be an
$M$-coalgebra, $a \in \scat{A}$, and $x \in Xa$.  A \demph{resolution} of
$x$ is a diagram
\begin{equation}
\label{eq:half-resolution}
\cdots \gobymod{m_2} a_1 \gobymod{m_1} a_0 = a
\end{equation}
together with a sequence $(x_n)_{n\in\nat}$ such that $x_n \in Xa_n$, $x_0
= x$, and
\[
\xi(x_n) = m_{n+1} \otimes x_{n+1}
\]
for all $n\in\nat$.  I will also call $(x_n)_{n\in\nat}$ a
\demph{resolution of $x$ along} the diagram~\bref{eq:half-resolution}.
Clearly every element $x$ of a coalgebra has at least one resolution.

\begin{lemma}[Direction of resolution]
\label{lemma:uniquedirectionofresolution}
Let $(X, \xi)$ be a nondegenerate $M$-coalgebra, $a \in \scat{A}$, and $x
\in Xa$.  If $x$ has resolutions along both
\[
( \cdots \gobymod{m_2} a_1 \gobymod{m_1} a_0 = a )
\diagspace
\textrm{and}
\diagspace
( \cdots \gobymod{m'_2} a'_1 \gobymod{m'_1} a'_0 = a )
\]
then
\[
\bigcompt{ \cdots \gobymod{m_2} a_1 \gobymod{m_1} a_0 }
=
\bigcompt{ \cdots \gobymod{m'_2} a'_1 \gobymod{m'_1} a'_0 }.
\]
\end{lemma}

\begin{proof}
Choose sequences $(x_n)_{n\in\nat}$ and $(x'_n)_{n\in\nat}$ resolving $x$
along the two diagrams respectively.  I construct by induction a
commutative diagram
\[
\begin{diagram}
\cdots		&\rMod^{m_3}	&a_2		&\rMod^{m_2}	&
a_1		&\rMod^{m_1}	&a_0 = a	\\
		&		&\uTo>{f_2}	&		&
\uTo>{f_1}	&		&\uTo>{f_0 = 1_a}\\
\cdots		&\rMod^{p_3}	&b_2		&\rMod^{p_2}	&
b_1		&\rMod^{p_1}	&b_0 = a	\\
		&		&\dTo>{f'_2}	&		&
\dTo>{f'_1}	&		&\dTo>{f'_0 = 1_a}\\
\cdots		&\rMod_{m'_3}	&a'_2		&\rMod_{m'_2}	&
a'_1		&\rMod_{m'_1}	&a'_0 = a	\\
\end{diagram}
\]
and a sequence $(y_n)_{n\in\nat}$ such that $y_n \in Xb_n$, $f_n y_n =
x_n$, and $f'_n y_n = x'_n$ for each $n\in\nat$.  For the base step, take
$y_0 = x$.  For the inductive step, let $n\in\nat$ and suppose that $b_n$,
$f_n$, $f'_n$, and $y_n$ are constructed.  We may write
\[
\xi(y_n) = (c \gobymod{q} b_n) \otimes z
\]
with $z \in Xc$.  Then 
\[
m_{n+1} \otimes x_{n+1}	
=
\xi(x_n)
=	
\xi(f_n y_n)
=
f_n \xi(y_n)
=
(f_n q) \otimes z,  
\]
so by nondegeneracy of $X$ and Lemma~\ref{lemma:equalityinMtensorX},
there exist a commutative diagram as labelled~(a) below and an element $w
\in Xd$ such that $gw = z$ and $hw = x_{n+1}$:
\[
\begin{diagram}
	&	&a_{n+1}&	&\rMod^{m_{n+1}}&&a_n	\\
	&	&\uTo<h	&	&	&	&	\\
	&	&d	&	&\textrm{(a)}&	&\uTo>{f_n}\\
	&\ruTo<k&	&\rdTo>g&	&	&	\\
b_{n+1}	&	&\textrm{(b)}&	&c	&\rMod^q&b_n	\\
	&\rdTo<{k'}&	&\ruTo>{g'}&	&	&	\\
	&	&d'	&	&\textrm{(a$'$)}&&\dTo>{f'_n}\\
	&	&\dTo<{h'}&	&	&	&	\\
	&	&a'_{n+1}&	&\rMod_{m'_{n+1}}&&a'_n.\\
\end{diagram}
\]
Similarly, there exist a commutative diagram~(a$'$) and $w' \in Xd'$ such
that $g'w' = z$ and $h'w' = x'_{n+1}$.  So by nondegeneracy of $X$
(condition \cstyle{ND1}), there exist a commutative square~(b) and $y_{n+1}
\in Xb_{n+1}$ such that $k y_{n+1} = w$ and $k' y_{n+1} = w'$.  Putting
$p_{n+1} = qgk$, $f_{n+1} = hk$, and $f'_{n+1} = h'k'$ completes the
induction.  \done
\end{proof}

Now consider resolutions in the coalgebra $(I, \iota)$.  Given any $a \in
\scat{A}$ and
\begin{equation}
\label{eq:diag-to-resolve}
( \cdots \gobymod{m_2} a_1 \gobymod{m_1} a_0 )
\in
\catI a,
\end{equation}
there is a canonical%
\label{p:tautological}
resolution of
$
\compt{
\cdots \gobymod{m_2} a_1 \gobymod{m_1} a_0
}
\in
Ia,
$
consisting of~\bref{eq:diag-to-resolve} itself together with 
$
\compt{
\cdots \gobymod{m_{n+2}} a_{n+1} \gobymod{m_{n+1}} a_n
}
\in
Ia_n
$
as `$x_n$'.

\begin{propn}[Double complex]
\label{propn:doublecomplex}
Let
\[
\begin{diagram}
\ddots	&		&\vdots	&		&
	&		&\vdots		\\
	&		&	&		&
	&		&\dMod>{m_3}	\\
\cdots	&\rMod^{m_2^3}	&a_2^2	&\rMod^{m_2^2}	&
a_2^1	&\rMod^{m_2^1}	&a_2^0		\\
	&		&	&		&
	&		&\dMod>{m_2}	\\
\cdots	&\rMod^{m_1^3}	&a_1^2	&\rMod^{m_1^2}	&
a_1^1	&\rMod^{m_1^1}	&a_1^0		\\
	&		&	&		&
	&		&\dMod>{m_1}	\\
\cdots	&\rMod^{m_0^3}	&a_0^2	&\rMod^{m_0^2}	&
a_0^1	&\rMod^{m_0^1}	&a_0^0		\\
\end{diagram}
\]
be a diagram satisfying
\[
\bigcompt{
\cdots 
\gobymod{m_n^3} a_n^2
\gobymod{m_n^2} a_n^1
\gobymod{m_n^1} a_n^0
}
=
\bigcompt{
\cdots 
\gobymod{m_{n+1}^2} a_{n+1}^1
\gobymod{m_{n+1}^1} a_{n+1}^0
\gobymod{m_{n+1}} a_n^0
}
\]
for all $n\in\nat$.  Then
\begin{equation}
\label{eq:dbl-cx}
\bigcompt{
\cdots \gobymod{m_0^3} a_0^2 \gobymod{m_0^2} a_0^1 \gobymod{m_0^1} a_0^0
}
=
\bigcompt{
\cdots \gobymod{m_3} a_2^0 \gobymod{m_2} a_1^0 \gobymod{m_1} a_0^0
}.
\end{equation}
\end{propn}

\begin{proof}
The left-hand side of~\bref{eq:dbl-cx} can be resolved canonically along
\[
\cdots \gobymod{m_0^2} a_0^1 \gobymod{m_0^1} a_0^0.
\]
It also has a resolution $(x_n)_{n\in\nat}$ along
\[
\cdots \gobymod{m_2} a_1^0 \gobymod{m_1} a_0^0,
\]
where
\[
x_n 
=
\bigcompt{
\cdots \gobymod{m_n^2} a_n^1 \gobymod{m_n^1} a_n^0
}
\in 
Ia_n^0,
\]
since by hypothesis
\begin{eqnarray*}
\iota x_n	&
=		&
\iota
\bigcompt{
\cdots \gobymod{m_{n+1}^1} a_{n+1}^0 \gobymod{m_{n+1}} a_n^0
}		\\
		&
=		&
m_{n+1} \otimes x_{n+1}.
\end{eqnarray*}
The result follows from nondegeneracy of $I$ and
Lemma~\ref{lemma:uniquedirectionofresolution}.  \done
\end{proof}

\begin{thm}[Universal solution in $\Set$]
\label{thm:universalsolutioninSet}
$(I, \iota)$ is the universal solution of $(\scat{A}, M)$ in $\Set$.
\end{thm}

\begin{proof}
Let $(X, \xi)$ be a nondegenerate coalgebra.  We have to show that there is
a unique map $(X, \xi) \go (I, \iota)$.

\paragraph*{Existence} Given any $a \in \scat{A}$ and $x \in Xa$, we may
choose a resolution
\begin{equation}
\label{eq:resolution}
( 
( \cdots \gobymod{m_2} a_1 \gobymod{m_1} a_0 ),
( x_0, x_1, x_2, \ldots )
)  
\end{equation}
of $x$ and put
\[
\ovln{\xi}_a(x) 
=
\bigcompt{
\cdots \gobymod{m_2} a_1 \gobymod{m_1} a
}
\in 
Ia.
\]
This defines for each $a$ a function $\ovln{\xi}_a: Xa \go Ia$, which by
Lemma~\ref{lemma:uniquedirectionofresolution} is independent of choice of
resolution.  I claim that $\ovln{\xi}$ is a map $(X, \xi) \go (I, \iota)$
of coalgebras.  First, it is a natural transformation, that is, if $a
\goby{f} a'$ is a map in $\scat{A}$ and $x \in Xa$ then $\ovln{\xi}_{a'}
(fx) = f \ovln{\xi}_a (x)$.  For choose a resolution~\bref{eq:resolution}
of $x$: then
\[
( 
( \cdots \gobymod{m_2} a_1 \gobymod{f m_1} a' ),
( fx, x_1, x_2, \ldots )
)  
\]
is a resolution of $fx$, so
\begin{eqnarray*}
\ovln{\xi}_{a'} (fx)	&
=	&
\bigcompt{ \cdots \gobymod{m_2} a_1 \gobymod{f m_1} a' }\\
	&
=	&
f \bigcompt{ \cdots \gobymod{m_2} a_1 \gobymod{m_1} a }\\
	&
=	&
f \ovln{\xi}_a (x).
\end{eqnarray*}
Second, $\ovln{\xi}$ is a map of coalgebras, that is, if $a \in \scat{A}$
and $x \in Xa$ then
\[
(M \otimes \ovln{\xi})_a \xi_a (x) 
= 
\iota_a \ovln{\xi}_a (x).
\]
For choose a resolution~\bref{eq:resolution} of $x$: then
\[
( 
( \cdots \gobymod{m_3} a_2 \gobymod{m_2} a_1 ),
( x_1, x_2, x_3, \ldots )
)  
\]
is a resolution of $x_1$, so
\begin{eqnarray*}
(M \otimes \ovln{\xi})_a \xi_a (x)	&
=	&
(M \otimes \ovln{\xi})_a (m_1 \otimes x_1)			\\
	&
=	&
m_1 \otimes \ovln{\xi}_{a_1} (x_1)				\\
	&
=	&
m_1 \otimes \bigcompt{ \cdots \gobymod{m_2} a_1 }		\\
	&
=	&
\iota_a \bigcompt{ \cdots \gobymod{m_2} a_1 \gobymod{m_1} a }	\\
	&
=	&
\iota_a \ovln{\xi}_a (x).
\end{eqnarray*}

\paragraph*{Uniqueness} Let $\twid{\xi}: (X, \xi) \go (I, \iota)$ be a map of
coalgebras, $a \in \scat{A}$, and $x \in Xa$.  We show that $\twid{\xi}_a
(x) = \ovln{\xi}_a (x)$.

Choose a resolution~\bref{eq:resolution} of $x$, and for each $n \in \nat$,
write
\[
\twid{\xi}_{a_n} (x_n)
=
\bigcompt{
\cdots \gobymod{m_n^2} a_n^1 \gobymod{m_n^1} a_n^0 = a_n
}.
\]
For each $n\in\nat$, we have
\[
(M \otimes \twid{\xi})_{a_n} \xi_{a_n} (x_n)
=
\iota_{a_n} \twid{\xi}_{a_n} (x_n)
\]
by definition of map of coalgebras.  On the other hand,
\begin{eqnarray*}
(M \otimes \twid{\xi})_{a_n} \xi_{a_n} (x_n)	&
=	&
(M \otimes \twid{\xi})_{a_n} (m_{n+1} \otimes x_{n+1})	\\
	&
=	&
m_{n+1} \otimes \twid{\xi}_{a_{n+1}} (x_{n+1})		\\
	&
=	&
m_{n+1} \otimes 
\bigcompt{ 
\cdots 
\gobymod{m_{n+1}^2} a_{n+1}^1 \gobymod{m_{n+1}^1} a_{n+1}^0 = a_{n+1}
}	\\
	&
=	&
\iota_{a_n}
\bigcompt{
\cdots 
\gobymod{m_{n+1}^1} a_{n+1}^0 = a_{n+1} \gobymod{m_{n+1}} a_n
}.
\end{eqnarray*}
Since $\iota_{a_n}$ is injective (Proposition~\ref{propn:Iisafixedpoint}),
\begin{eqnarray*}
\bigcompt{
\cdots 
\gobymod{m_{n+1}^1} a_{n+1}^0 = a_{n+1} \gobymod{m_{n+1}} a_n
}	&
=	&
\twid{\xi}_{a_n} (x_n)	\\
	&
=	&
\bigcompt{
\cdots \gobymod{m_n^2} a_n^1 \gobymod{m_n^1} a_n^0 = a_n
}
\end{eqnarray*}
for each $n\in\nat$.  So Proposition~\ref{propn:doublecomplex} applies, and
\[
\bigcompt{
\cdots \gobymod{m_0^2} a_0^1 \gobymod{m_0^1} a_0^0 
}
=
\bigcompt{
\cdots \gobymod{m_2} a_1 \gobymod{m_1} a_0
}
\in
Ia,
\]
that is, $\twid{\xi}_a(x) = \ovln{\xi}_a (x)$, as required.
\done
\end{proof}

\section{Topological proofs}
\label{sec:Top-proofs}

Fix a self-similarity system $(\scat{A}, M)$.  In this section I show that
if $(\scat{A}, M)$ satisfies $\So$ then $(I, \iota)$, with the topology
defined in~\S\ref{sec:univ-soln}, is an $M$-coalgebra in $\Top$, and indeed
the universal solution in $\Top$.

The proof involves an analysis of equality in $Ia$, that is, of the
possible representations in $\catI a$ of a given element of $Ia$.

\begin{lemma}
\label{lemma:Inafinite}
For each $a\in\scat{A}$ and $n\in\nat$, the category $\catI_n a$ is
finite. 
\end{lemma}

\begin{proof}
Follows from finiteness of $M$.
\done
\end{proof}

The last paragraph of~\S\ref{sec:univ-soln} shows that there may be objects
$\tau$, $\tau'$ of $\catI a$ such that $\pr_n (\tau)$ and $\pr_n (\tau')$
are in the same connected-component of $\catI_n a$ for all $n$ but $\tau$
and $\tau'$ are in different connected-components of $\catI a$.  However,
$\tau$ and $\tau'$ are in the same component if $\pr_n(\tau)$ and
$\pr_n(\tau')$ can be connected by a diagram \emph{of the same shape for
each $n$}:
\begin{lemma}[Equality in $Ia$]
\label{lemma:equalityinIa}
Let $a \in \scat{A}$ and $\tau, \tau' \in \catI a$.  Then $\compt{\tau} =
\compt{\tau'} \in Ia$ if and only if there exist a finite connected
category $\scat{K}$, objects $k, k' \in \scat{K}$, and for each $n \geq
1$, a functor $F_n: \scat{K} \go \catI_n a$ such that $F_n k = \pr_n
\tau$ and $F_n k' = \pr_n \tau'$.
\end{lemma}

\begin{proof}
`Only if' is simple.  For `if' we use Lemma~\ref{lemma:factorization}.
Take $\scat{J}$ to be the discrete category on two objects $j$, $j'$,
define $P: \scat{J} \go \scat{K}$ by $P(j) = k$ and $P(j') = k'$, take
$\cat{L}$ to be the limit $\catI a$ of the categories $\cat{L}_n = \catI_n
a$, and take $D(j) = \tau$ and $D(j') = \tau'$.  Then for each $n\geq 1$,
the set of factorizations in Lemma~\ref{lemma:factorization} is nonempty by
hypothesis, and finite since $\scat{K}$ and $\catI_n a$ are finite, so
there is a functor $G: \scat{K} \go \catI a$ such that $Gk = \tau$ and $Gk'
= \tau'$.  But $\scat{K}$ is connected, so $\compt{\tau} = \compt{\tau'}$.
\done
\end{proof}

Assume from now on that $(\scat{A}, M)$ satisfies \So.

Each module element $q: d \gomod c$ induces a function $\phi_q: Id \go Ic$
by
\[
\phi_q
\bigcompt{
\cdots \gobymod{r_2} d_1 \gobymod{r_1} d
}
=
\bigcompt{
\cdots \gobymod{r_2} d_1 \gobymod{r_1} d \gobymod{q} c
},
\]
or equivalently,
\[
\label{p:phi}
\phi_q
=
\left(
Id \goby{q\otimes\dashbk} 
(M \otimes I)c \rTo^{\iota^{-1}_c}_{\diso}
Ic
\right)
\]
where the first map is the coprojection
\[
\mr{copr}_q 
=
q \otimes \dashbk:
Id 
\go
\int^{d'} M(d', c) \times Id' 
=
(M \otimes I) c.
\]
For each $a \in \scat{A}$, the topology on $Ia$ is generated%
\label{p:topo-gen}
by taking $\phi_{p_1} \phi_{p_2} \cdots \phi_{p_n} (Ib_n)$ to be a closed
subset of $Ia$ whenever $n\in\nat$ and
\[
(b_n \gobymod{p_n} \cdots \gobymod{p_2} b_1 \gobymod{p_1} b_0)
\in
\catI_n a.
\]
It is shown in the following pages that $(I, \iota)$ is an $M$-coalgebra in
$\Top$, and along the way that each $Ia$ is compact Hausdorff.

We start with the Hausdorff property.  Define, for each $n\in\nat$ and
$a\in\scat{A}$, a binary relation $R_n^a$ on $Ia$ by
\[
R_n^a 
=
\bigcup
\{
(\phi_{p_1} \cdots \phi_{p_n} (Ib_n))^2
\such
(b_n \gobymod{p_n} \cdots \gobymod{p_1} b_0)
\in 
\catI_n a
\}
\sub
Ia \times Ia.
\]
Equivalently, $(t, t') \in R_n^a$ when there exists $(b_n \gobymod{p_n}
\cdots \gobymod{p_1} b_0)$ such that $t$ and $t'$ can both be written in
the form 
\[
\bigcompt{
\cdots \gomod \cdot \gomod 
b_n \gobymod{p_n} \cdots \gobymod{p_1} b_0
}.
\]
As a subset of $Ia \times Ia$, $R_n^a$ is closed (by finiteness of
$\catI_n a$).  As a relation, $R_n^a$ is reflexive and symmetric, but not
in general transitive: for instance, in the Freyd self-similarity system,
$
R_1^1 
=
[0, \half]^2 \cup [\half, 1]^2
% [0, \dhalf]^2 \cup [\dhalf, 1]^2
\sub
[0, 1]^2.
$

For any set $S$, write $\Delta_S = \{ (s, s) \such s \in S \} \sub S \times
S$. 

\begin{lemma}[Relations determine equality]
\label{lemma:relationsdetermineequality}
$\bigcap_{n\in\nat} R_n^a = \Delta_{Ia}$ for each $a$.
\end{lemma}

\begin{proof}
Certainly $\bigcap_{n\in\nat} R_n^a \supseteq \Delta_{Ia}$.  Conversely, let
$(t, t') \in \bigcap_{n\in\nat} R_n^a$, writing
\begin{eqnarray*}
t	&=	&
\compt{
\cdots \gobymod{m_2} a_1 \gobymod{m_1} a_0 = a
},	\\
t'	&=	&
\compt{
\cdots \gobymod{m'_2} a'_1 \gobymod{m'_1} a'_0 = a
}.
\end{eqnarray*}
For each $n\in\nat$, we may choose $(b_n \gobymod{p_n} \cdots \gobymod{p_1}
b_1) \in \catI_n a$ such that $t, t' \in \phi_{p_1} \cdots \phi_{p_n}
(Ib_n)$.  Since $\catI a$ is componentwise cofiltered, there is for each
$n\in\nat$ a span in $\catI a$ of the form
\[
\begin{diagram}
\cdots	&
\rMod^{m_{n+2}}	&a_{n+1}	&
\rMod^{m_{n+1}}	&a_n		&
\rMod^{m_n}	&\ 		&
\cdots	&
\		&\rMod^{m_1}	&a_0 = a	\\
	&
		&\uTo		&
		&\uTo		&
		&		&
	&
		&		&\uTo>{1_a}	\\
\cdots	&
\rMod		&\cdot		&
\rMod		&\cdot		&
\rMod		&\ 		&
\cdots	&
\		&\rMod		&a		\\
	&
		&\dTo		&
		&\dTo		&
		&		&
	&
		&		&\dTo>{1_a}	\\
\cdots	&
\rMod		&\cdot		&
\rMod		&b_n		&
\rMod_{p_n}	&\ 		&
\cdots	&
\		&\rMod_{p_1}	&b_0 = a,	\\
\end{diagram}
\]
hence a span in $\catI_n a$ of the form
\[
\begin{diagdiag}
	&	&\cdot	&	&	\\
	&\ldTo	&	&\rdTo	&	\\
\begin{diagram}
(a_n &\gobymod{m_n}& \cdots &\gobymod{m_1}& a_0)
\end{diagram}	&
		&	&	&
\begin{diagram}
(b_n &\gobymod{p_n} &\cdots &\gobymod{p_1}& b_0).
\end{diagram}	\\
\end{diagdiag}
\]
The same is true for $t'$, so for each $n\in\nat$ there is a diagram of the
form
\[
\begin{diagdiag}
	&	&\cdot	&	&	&	&\cdot	&	&	\\
	&\ldTo	&	&\rdTo	&	&\ldTo	&	&\rdTo	&	\\
\begin{diagram}
(a_n &\gobymod{m_n} &\cdots &\gobymod{m_1} &a_0)	
\end{diagram}	&
		&	&	&\cdot	&	&	&	&
\begin{diagram}
(a'_n &\gobymod{m'_n} &\cdots &\gobymod{m'_1}& a'_0)	
\end{diagram}\\
\end{diagdiag}
\]
in $\catI_n a$.  So by Lemma~\ref{lemma:equalityinIa}, taking $\scat{K}$
to be the evident category with four non-identity arrows, $t = t'$.
\done
\end{proof}

\begin{propn}[$Ia$ Hausdorff]
\label{propn:Iahausdorff}
$Ia$ is Hausdorff for all $a \in \scat{A}$. 
\end{propn}

\begin{proof}
$\Delta_{Ia}$ is the intersection of the closed subsets $R_n^a$ of $Ia
\times Ia$.
\done
\end{proof}

The next step is to consider $I: \scat{A} \go \Set$ as a quotient of
$\oba\catI$, the composite of $\catI: \scat{A} \go \Cat$ with the objects
functor $\ob: \Cat \go \Set$.  This functor $\oba\catI$ is nondegenerate
(even if \So\ does not hold), since
\[
\oba\catI  
\iso
\sum_b 
\ob(\catI b) \times M(b, \dashbk)
\]
and the class of nondegenerate functors is closed under sums.  Moreover,
$\oba\catI$ carries an $M$-coalgebra structure $\iota$, given by the usual
formula 
\[
\iota_{a_0}
( \cdots \gobymod{m_2} a_1 \gobymod{m_1} a_0 )
\ 
=
\ 
m_1
\otimes
(\cdots \gobymod{m_2} a_1)
\]
or equivalently by taking $\iota_a$ to be the composite
\begin{equation}	\label{eq:iota-as-comp}
\ob(\catI a)
\goiso
\sum_b M(b, a) \times \ob(\catI b)
\goby{\mr{canonical}}
\int^b
M(b, a) \times \ob(\catI b).
\end{equation}
So $(\oba\catI, \iota)$ is a coalgebra in $\Set$, and there is a canonical
map of coalgebras $\pi: (\oba\catI, \iota) \go (I, \iota)$.  

This construction can be topologized.  For each $a \in \scat{A}$, the set
$\ob(\catI a)$ is the limit of the diagram
\[
\cdots \go \ob(\catI_2 a) \go \ob(\catI_1 a),
\]
and giving each $\ob(\catI_n a)$ the discrete topology induces a topology
on $\ob(\catI a)$; in this way, $\oba\catI$ becomes a functor $\scat{A} \go
\Top$.  Each $\catI_n a$ is finite, so each $\ob(\catI a)$ is compact
Hausdorff, so $\oba\catI: \scat{A} \go \Top$ is nondegenerate.  Moreover,
the maps $\iota_a$ are continuous, since in~\bref{eq:iota-as-comp} the
first map is a homeomorphism and the second is a quotient map.  So
$(\oba\catI, \iota)$ is a coalgebra in $\Top$.  We will see that the maps
$\pi_a: \ob(\catI a) \go Ia$ exhibit each $Ia$ as not merely a
set-theoretic quotient of $\ob(\catI a)$, but a topological quotient.

For each $a \in \scat{A}$, the space $\ob(\catI a)$ is compact.  It is also
\demph{second countable} (has a countable basis of open sets).  Hence:
\begin{lemma}
\label{lemma:sequentialcompactness}
Every sequence in $\ob(\catI a)$ has a convergent subsequence.
\done 
\end{lemma}

The following result must hold if the spaces $Ia$ are to be compact.

\begin{lemma}[Intersections of basic closed sets]
\label{lemma:basicclosedsets}
Let $a \in \scat{A}$, $n \in \nat$,
\[
(\cdots \gobymod{m_2} a_1 \gobymod{m_1} a_0 )
\in
\catI a,
\diagspace
(b_n \gobymod{p_n} \cdots \gobymod{p_1} b_0)
\in
\catI_n a. 
\]
Then 
\[
\bigcompt{ \cdots \gobymod{m_2} a_1 \gobymod{m_1} a_0 } 
\in 
\phi_{p_1} \cdots \phi_{p_n} (Ib_n)
\]
if and only if
\[
\forall r \in \nat, 
\ 
\phi_{m_1} \cdots \phi_{m_r} (I a_r)
\cap
\phi_{p_1} \cdots \phi_{p_n} (I b_n)
\neq
\emptyset.
\]
\end{lemma}

\begin{proof}
`Only if' is trivial.  For `if', write $\alpha = ( \cdots \gobymod{m_1} a_0
)$.  By hypothesis, we may choose for each $r\in\nat$ objects
\begin{eqnarray*}
\alpha_r	&=	&
(
\cdots
\gobymod{m_{r+2}^r} a_{r+1}^r \gobymod{m_{r+1}^r} a_r \gobymod{m_r}
\cdots
\gobymod{m_1} a_0
),	\\
\beta_r	&=	&
(
\cdots
\gobymod{p_{n+2}^r} b_{n+1}^r \gobymod{p_{n+1}^r} b_n \gobymod{p_n}
\cdots
\gobymod{p_1} b_0
) 
\end{eqnarray*}
of $\catI a$ such that $\compt{\alpha_r} = \compt{\beta_r}$.  By
Proposition~\ref{propn:Iacomponentwisecofiltered}, there is for each
$r\in\nat$ a span
\[
\alpha_r 
\og \cdot \go 
\beta_r
\]
in $\catI a$.  By Lemma~\ref{lemma:sequentialcompactness},
$(\beta_r)_{r\in\nat}$ has a subsequence $(\beta_{r_k})_{k\in\nat}$ convergent
to $\beta$, say, and $\beta$ is of the form
\[
\cdots 
\gobymod{p_{n+1}} b_n \gobymod{p_n}
\cdots
\gobymod{p_1} b_0.
\]
For each $r\in\nat$, we may choose $k\in\nat$ such that $r_k \geq r$ and
$\pr_r (\beta_{r_k}) = \pr_r (\beta)$; we then have a span in $\catI a$ of
the form
\[
\alpha_{r_k}
\og \cdot \go
\beta_{r_k},
\]
hence, applying $\pr_r$, a span in $\catI_r a$ of the form
\[
\pr_r (\alpha)
\og \cdot \go
\pr_r (\beta).
\]
So by Lemma~\ref{lemma:equalityinIa}, $\compt{\alpha} = \compt{\beta} \in
\phi_{p_1} \cdots \phi_{p_n} (Ib_n)$.  \done
\end{proof}

\begin{propn}[$Ia$ as a quotient]
\label{propn:Iaasaquotient}
For each $a\in\scat{A}$, the canonical surjection $\pi_a: \ob(\catI a)
\go Ia$ is a topological quotient map.
\end{propn}

\begin{proof}
First, $\pi_a$ is continuous.  Let $n\in\nat$ and $(b_n \gobymod{p_n}
\cdots \gobymod{p_1} b_0) \in \catI_n a$; we must show that $\pi_a^{-1}
(\phi_{p_1} \cdots \phi_{p_n} (Ib_n))$ is a closed subset of $\ob(\catI
a)$.  By Lemma~\ref{lemma:basicclosedsets}, an element $( \cdots
\gobymod{m_1} a_0)$ of $\ob(\catI a)$ belongs to this subset if and only if
\begin{equation} \label{eq:nonempty-intersection}
\phi_{m_1} \cdots \phi_{m_r} (I a_r)
\cap
\phi_{p_1} \cdots \phi_{p_n} (I b_n)
\neq
\emptyset
\end{equation}
for all $r\in\nat$.  In other words,
\[
\pi_a^{-1} (\phi_{p_1} \cdots \phi_{p_n} (I b_n)) 
=
\bigcap_{r\in\nat} \pr_r^{-1} W_r
\]
where $\pr_r: \catI a \go \catI_r a$ and $W_r \sub \ob(\catI_r a)$ is the
set of elements $(a_r \gobymod{m_r} \cdots \gobymod{m_1} a_0)$
satisfying~\bref{eq:nonempty-intersection}.  But each $\ob(\catI_r a)$ is
discrete and each $\pr_r$ continuous, so $\bigcap_{r\in\nat} \pr_r^{-1}
W_r$ is closed, as required.

Second, $\pi_a$ is closed, since $\ob(\catI a)$ is compact and $Ia$ is
Hausdorff; and any continuous closed surjection is a quotient map.
\done
\end{proof}

\begin{cor}[$Ia$ compact]
\label{cor:Iacompact}
$Ia$ is compact for all $a \in \scat{A}$.  \done
\end{cor}

\begin{cor}[Properties of $If$]
\label{cor:propertiesofIf}
For each map $a \goby{f} a'$ in $\scat{A}$, the map $Ia \goby{If} Ia'$ is
continuous and closed.
\end{cor}

\begin{proof}
There is a commutative square
\[
\begin{diagram}
\ob(\catI a)	&\rQt^{\pi_a}	&Ia		\\
\dTo<{\ob(\catI f)}	&		&\dTo>{If}	\\
\ob(\catI a')	&\rQt_{\pi_{a'}}&Ia'		\\
\end{diagram}
\]
in which $\pi_a$ is a topological quotient map and $\ob(\catI f)$ and
$\pi_{a'}$ are continuous, so $If$ is also continuous.  $Ia$ is compact and
$Ia'$ Hausdorff, so $If$ is closed.  
\done
\end{proof}

\begin{cor}[$\iota_a$ continuous]
\label{cor:iotaacontinuous}
For each $a \in \scat{A}$, the map $\iota_a: Ia \go (M \otimes I) a$ is
continuous. 
\end{cor}

\begin{proof}
As for the previous lemma, using the square 
\[
\begin{diagram}
\ob(\catI a)		&\rQt^{\pi_a}			&Ia		\\
\dTo<{\iota_a}		&				&\dTo>{\iota_a}	\\
(M\otimes \oba\catI)a	&\rTo_{(M \otimes \pi)_a}	&(M \otimes I) a.\\
\end{diagram}
\]
\done
\end{proof}

\begin{cor}
\label{cor:topologicalcoalgebra}
$(I, \iota)$ is an $M$-coalgebra in $\Top$.  \done
\end{cor}

Our final task is to prove that for any $M$-coalgebra $(X, \xi)$ in
$\Top$, the unique map $\ovln{\xi}: (X, \xi) \go (I, \iota)$ of coalgebras
in $\Set$ is continuous.  To do this we show that the inverse image
$\ovln{\xi}_a^{-1} ( \phi_{p_1} \cdots \phi_{p_n} (I b_n) )$ of each basic
closed set is closed.  This inverse image is larger than it might appear.
For, given $x \in Xa$, one might imagine that if $\ovln{\xi}_a(x) \in
\phi_{p_1} \cdots \phi_{p_n} (I b_n)$ then $x$ can be resolved along some
diagram of the form
\[
\cdots \gomod \cdot \gomod b_n \gobymod{p_n} \cdots \gobymod{p_1} b_0;
\]
but while this is certainly a sufficient condition, it is not necessary:

\begin{example}
Let $(\scat{A}, M)$ be the Freyd self-similarity system.  Choose an
endpoint-preserving continuous map $\xi_1: [0, 1] \go [0, 2]$ satisfying
$\xi_1(2/3) = 2/3$; this gives an $M$-coalgebra structure $\xi$ on $X = (
\{ \star \} \parpair{0}{1} [0, 1] )$.  Now $2/3 \in X1$ has a resolution
along
\[
\cdots \gobymod{[0, \half]} 1 \gobymod{[0, \half]} 1, 
\]
so $\ovln{\xi}_1 (2/3) = 0$, so $\ovln{\xi}_1 (2/3) \in \phi_{m_1} (I0)$
where $m_1 = 0: 0 \gomod 1$.  But $2/3$ has no resolution ending in
$m_1$.
\end{example}

To describe this inverse image we need some notation.  Given an
$M$-coalgebra $X = (X, \xi)$ in $\Top$, $r \in \nat$, and $(a_r
\gobymod{m_r} \cdots \gobymod{m_1} a_0) \in \catI_r a$, let
\[
\Res_X (a_r \gobymod{m_r} \cdots \gobymod{m_1} a_0) 
\sub
Xa
\]
be the set of all $x \in Xa$ having a resolution along some diagram of the
form
\[
\cdots \gomod \cdot \gomod 
a_r \gobymod{m_r} \cdots \gobymod{m_1} a_0.
\]

\begin{lemma}[Inverse image]
\label{lemma:lordy}
Let $a \in \scat{A}$, $n \in \nat$, and $(b_n \gobymod{p_n} \cdots
\gobymod{p_1} b_0) \in \catI_n a$; let $(X, \xi)$ be an $M$-coalgebra in
$\Set$ and $x \in Xa$.  Then $x \in \ovln{\xi}_a^{-1} ( \phi_{p_1} \cdots
\phi_{p_n} (I b_n) )$ if and only if
\begin{condition}
for all $r \in \nat$, there exists $(a_r \gobymod{m_r} \cdots \gobymod{m_1}
a_0) \in \catI_r a$ such that $x \in \Res_X (a_r \gobymod{m_r} \cdots
\gobymod{m_1} a_0)$ and
$
\phi_{m_1} \cdots \phi_{m_r} (I a_r) 
\cap
\phi_{p_1} \cdots \phi_{p_n} (I b_n)
\neq
\emptyset.
$
\end{condition}
\end{lemma}

\begin{proof}
For `only if', choose $( \cdots \gobymod{m_1} a_0) \in \catI a$
along which $x$ can be resolved.  Then
\[
\bigcompt{ \cdots \gobymod{m_1} a_0 }
=
\ovln{\xi}_a (x)
\in
\phi_{p_1} \cdots \phi_{p_n} (I b_n)
\]
and the result follows.

For `if', choose for each $r \in \nat$ an element $( \cdots \gobymod{m^r_1}
a^r_0) \in \catI a$ along which $x$ can be resolved and such that
\[
\phi_{m^r_1} \cdots \phi_{m^r_r} (I a^r_r) 
\cap
\phi_{p_1} \cdots \phi_{p_n} (I b_n)
\neq
\emptyset.
\]
By compactness of $Ia$, 
\begin{equation} \label{eq:closed-intersection}
\bigcap_{r \in \nat}
\phi_{m^r_1} \cdots \phi_{m^r_r} (I a^r_r) 
\cap
\phi_{p_1} \cdots \phi_{p_n} (I b_n)
\end{equation}
is nonempty.  But $\bigcap_{r \in \nat} \phi_{m^r_1} \cdots \phi_{m^r_r} (I
a^r_r)$ has cardinality at most $1$ by
Lemma~\ref{lemma:relationsdetermineequality}, and contains $\ovln{\xi}_a
(x)$ since $x$ can be resolved along $( \cdots \gobymod{m^r_1} a^r_0)$ for
each $r$, so the set~\bref{eq:closed-intersection} is $\{ \ovln{\xi}_a (x)
\}$.  The result follows.  \done
\end{proof}

This says what the inverse images of the basic closed sets are, and we now
prepare to show that they are closed.
\begin{lemma}[Coprojections closed]
\label{lemma:coprojectionsclosed}
Let $X: \scat{A} \go \Top$ be a nondegenerate functor and $m: b \gomod a$ a
module element.  Then the coprojection
\[
m \otimes \dashbk:
Xb \go (M \otimes X)a
\]
is closed.
\end{lemma}

\begin{proof}
Trivially, $m \otimes \dashbk$ is the composite
\[
Xb
\goby{(m, \dashbk)}
\sum_{b'} M(b', a) \times Xb'
\rQt^{\mr{canonical}}
\int^{b'} M(b', a) \times Xb'
=
(M \otimes X)a.
\]
The first map is closed since it is a coproduct-coprojection, and the
second is closed by Lemma~\ref{lemma:closedquotientmap}.  
\done
\end{proof}

\begin{lemma}
\label{lemma:resolutionsetsclosed}
Let $(X, \xi)$ be an $M$-coalgebra in $\Top$, $r\in\nat$, and $(a_r
\gobymod{m_r} \cdots \gobymod{m_1} a_0) \in \catI_r a$.  Then $\Res_X (a_r
\gobymod{m_r} \cdots \gobymod{m_1} a_0)$ is a closed subset of $Xa$.
\end{lemma}

\begin{proof}
When $r = 0$ this is trivial.  Suppose inductively that the result holds
for $r \in \nat$, and let $(a_{r+1} \gobymod{m_{r+1}} \cdots \gobymod{m_1}
a_0) \in \catI_{r+1} a$.  An element $x \in Xa$ is in $\Res_X(a_{r+1}
\gobymod{m_{r+1}} \cdots \gobymod{m_1} a_0)$ if and only if there exists
$x_1 \in Xa_1$ such that
\[
\xi_a (x)
=
m_1 \otimes x_1
\diagspace
\textrm{and}
\diagspace
x_1 \in \Res_X (a_{r+1} \gobymod{m_{r+1}} \cdots \gobymod{m_2} a_1), 
\]
so
\[
\Res_X(a_{r+1} \gobymod{m_{r+1}} \cdots \gobymod{m_1} a_0)
= 
\xi_a^{-1}
\left(
m_1 \otimes \Res_X (a_{r+1} \gobymod{m_{r+1}} \cdots \gobymod{m_2} a_1)
\right).
\]
But $\Res_X (a_{r+1} \gobymod{m_{r+1}} \cdots \gobymod{m_2} a_1)$ is closed
by inductive hypothesis, $m_1 \otimes \dashbk$ is closed by
Lemma~\ref{lemma:coprojectionsclosed}, and $\xi_a$ is continuous,
so the induction is complete.
\done
\end{proof}

\begin{thm}[Universal solution in $\Top$]
\label{thm:universalsolutioninTop}
$(I, \iota)$ is the universal solution of $(\scat{A}, M)$ in $\Top$.  
\end{thm}

\begin{proof}
It remains to show that for any $M$-coalgebra $(X, \xi)$ in $\Top$, the
unique map $\ovln{\xi}: (X, \xi) \go (I, \iota)$ is continuous in each
component.  So let $a \in \scat{A}$, $n \in \nat$, and $(b_n \gobymod{p_n}
\cdots \gobymod{p_1} b_0) \in \catI_n a$; we must show that
\[
\ovln{\xi}_a^{-1} (\phi_{p_1} \cdots \phi_{p_n} (I b_n))
\sub
Xa
\]
is closed.  For each $r \in \nat$, write
\[
V_r
=
\bigcup
\Res_X (a_r \gobymod{m_r} \cdots \gobymod{m_1} a_0)
\]
where the union is over all $(a_r \gobymod{m_r} \cdots \gobymod{m_1} a_0)
\in \catI_r a$ such that
\[
\phi_{m_1} \cdots \phi_{m_r} (I a_r)
\cap
\phi_{p_1} \cdots \phi_{p_n} (I b_n)
\neq 
\emptyset.
\]
Then Lemma~\ref{lemma:lordy} says that
$
\ovln{\xi}_a^{-1} (\phi_{p_1} \cdots \phi_{p_n} (I b_n))
=
\bigcap_{r \in \nat} V_r,
$
so it is enough to prove that each $V_r$ is closed, and this follows from
Lemmas~\ref{lemma:Inafinite} and~\ref{lemma:resolutionsetsclosed}.  \done
\end{proof}

%% file: solvability.tex
\section{Appendix: Solvability}
\label{app:solv}

Here we finish the proof of
\begin{thm}[Existence of universal solution]
\label{thm:existenceofuniversalsolution}
The following are equivalent for a self-similarity system $(\scat{A}, M)$:
\begin{enumerate}
\item	\label{item:ex-S}
$(\scat{A}, M)$ satisfies \So
\item	\label{item:ex-Top}
$(\scat{A}, M)$ has a universal solution in $\Top$
\item	\label{item:ex-Set}
$(\scat{A}, M)$ has a universal solution in $\Set$.
\end{enumerate}
\end{thm}
We proved \bref{item:ex-S}$\implies$\bref{item:ex-Top} in
\S\ref{sec:Top-proofs} and
\bref{item:ex-Top}$\implies$\bref{item:ex-Set} as
Proposition~\ref{propn:coalgebrasinTopandSet}, so it remains to prove
\bref{item:ex-Set}$\implies$\bref{item:ex-S}.

Fix a self-similarity system $(\scat{A}, M)$.  In this appendix,
`$M$-coalgebra' means `$M$-coalgebra in $\Set$'.

We begin by constructing representable-type coalgebras and proving a
Yoneda-type lemma.  Take $(a_\blb, m_\blb) = ( \cdots \gobymod{m_2} a_1
\gobymod{m_1} a_0 ) \in \elt{\catI}$.  Then there is a functor
\[
H^{(a_\sblb, m_\sblb)}
=
\sum_{n\in\nat}
\scat{A} (a_n, \dashbk):
\scat{A} \go \Set,
\]
and any representable functor is flat, so $H^{(a_\sblb, m_\sblb)}$ is
nondegenerate by Theorem~\ref{thm:componentwiseflatness}.  Also
\[
(M \otimes H^{(a_\sblb, m_\sblb)}) b
\iso
\sum_{n\in\nat} (M \otimes \scat{A}(a_n, \dashbk)) b
\iso	
\sum_{n\in\nat} M(a_n, b),
\]
so an $M$-coalgebra structure on $H^{(a_\sblb, m_\sblb)}$ amounts to a
natural transformation
\[
\sum_{n\in\nat} \scat{A}(a_n, \dashbk)
\go
\sum_{n\in\nat} M(a_n, \dashbk).
\]
There is a unique such transformation sending $1_{a_n}$ to $m_{n+1} \in
M(a_{n+1}, a_n)$ for each $n\in\nat$; let $\theta^{(a_\sblb, m_\sblb)}$
be the corresponding coalgebra structure on $H^{(a_\sblb, m_\sblb)}$.

This defines an $M$-coalgebra $(H^{(a_\sblb, m_\sblb)}, \theta^{(a_\sblb,
m_\sblb)})$ for each object $(a_\blb, m_\blb)$ of $\elt{\catI}$.  In
fact there is a functor
\[
(H^\blb, \theta^\blb):
\elt{\catI}^\op
\go
\Coalg{M}{\Set},
\]
since any map $f_\blb: (a_\blb, m_\blb) \go (a'_\blb, m'_\blb)$ in
$\elt{\catI}$ induces a map
\[
H^{(a'_\sblb, m'_\sblb)}
=
\sum_{n\in\nat} \scat{A}(a'_n, \dashbk)
\goby{\sum f_n^*}
\sum_{n\in\nat} \scat{A}(a_n, \dashbk)
=
H^{(a_\sblb, m_\sblb)}
\]
respecting the coalgebra structures. 

\begin{lemma}[`Yoneda']
\label{lemma:Yoneda}
There is a bijection
\begin{eqnarray}
	&
	&
\Coalg{M}{\Set}
\left(
(H^{(a_\sblb, m_\sblb)}, \theta^{(a_\sblb, m_\sblb)}),
(X, \xi)
\right)
\nonumber	\\
	&
\iso	&
\{
\textrm{\upshape resolutions along }
(a_\blb, m_\blb)
\textrm{\upshape\ in }
(X, \xi)
\}
\label{eq:Yoneda}
\end{eqnarray}
natural in $(a_\blb, m_\blb) \in \elt{\catI}$ and $(X, \xi) \in
\Coalg{M}{\Set}$.  More precisely, if $x \in Xa_0$ then the coalgebra maps
$H^{(a_\sblb, m_\sblb)} \go X$ sending $1_{a_0}$ to $x$ correspond to the
resolutions of $x$ along $(a_\blb, m_\blb)$.
\end{lemma}

\begin{proof}
A natural transformation $\alpha: H^{(a_\sblb, m_\sblb)} \go X$ amounts to a
sequence $(x_n)_{n\in\nat}$ with $x_n \in Xa_n$, by the standard Yoneda
Lemma.  It is a map of coalgebras if and only if
\[
\begin{diagram}[height=6ex]
\sum_{n\in\nat} \scat{A}(a_n, \dashbk)	&
\rTo^\alpha		&
X	\\
\dTo<{\theta^{(a_\sblb, m_\sblb)}}	&
		&
\dTo>\xi\\
\sum_{n\in\nat} M(a_n, \dashbk)		&
\rTo_{M \otimes \alpha}	&
M \otimes X\\
\end{diagram}
\]
commutes, if and only if this diagram commutes when we take $1_{a_n}$ at
the top-left corner for every $n\in\nat$.  We have
\[
\begin{diagram}
1_{a_n}		&\rGoesto	&x_n			\\
\dGoesto	&		&\dGoesto		\\
		&		&\xi(x_n)		\\
m_{n+1}		&\rGoesto	&m_{n+1} \otimes x_{n+1},	
\end{diagram}
\]
so $\alpha$ is a map of coalgebras just when $\xi(x_n) = m_{n+1} \otimes
x_{n+1}$ for all $n\in\nat$.  A coalgebra map $(H^{(a_\sblb, m_\sblb)},
\theta^{(a_\sblb, m_\sblb)}) \go (X, \xi)$ therefore amounts to a sequence
$(x_n)_{n\in\nat}$ satisfying $\xi(x_n) = m_{n+1} \otimes x_{n+1}$ for all
$n$, that is, a resolution along $(a_\blb, m_\blb)$ in $(X, \xi)$.  This
establishes the bijection~\bref{eq:Yoneda}; naturality follows from the
naturality in the standard Yoneda Lemma.  \done
\end{proof}

We have met one other $M$-coalgebra: $(\oba\catI, \iota)$, constructed
after Proposition~\ref{propn:Iahausdorff}.  (Recall
from~\S\ref{sec:Set-proofs} that $(I, \iota)$ is only known to be a
(nondegenerate) coalgebra if \So\ holds.)

\begin{cor}[Tautological map]
\label{cor:tautologicalmap}
For each $(a_\blb, m_\blb) \in \elt{\catI}$ there is a canonical map of
$M$-coalgebras  
\[
\kappa^{(a_\sblb, m_\sblb)}:
(H^{(a_\sblb, m_\sblb)}, \theta^{(a_\sblb, m_\sblb)})
\go
(\oba \catI, \iota)
\]
sending $1_{a_0}$ to $(a_\blb, m_\blb)$.
\end{cor}
\begin{proof}
Every $(a_\blb, m_\blb) \in \elt{\catI}$, regarded as an element of
$(\oba\catI) a_0$, has a tautological resolution in $\oba\catI$.  By
Lemma~\ref{lemma:Yoneda}, the corresponding map $\kappa^{(a_\sblb, m_\sblb)}$
of coalgebras sends $1_{a_0}$ to $(a_\blb, m_\blb)$.  \done
\end{proof}

We can now finish the proof of
Theorem~\ref{thm:existenceofuniversalsolution}.

Suppose that $(\scat{A}, M)$ has a universal solution $(J, \gamma)$.  Then
there is a map $\beta: (\oba\catI, \iota) \go (J, \gamma)$ of
$M$-coalgebras.  I claim that the natural transformation $\beta$ can be
factorized as
\[
\begin{diagram}
\oba \catI 		&		&
\rTo^\beta		&
			&J,		\\
			&\rdQt<\pi	&
			&
\ruGet>{\ovln{\beta}}	&		\\
			&		&
I			&
			&		\\
\end{diagram}
\]
or equivalently that for each $a \in \scat{A}$, the function $\beta_a:
\ob(\catI a) \go Ja$ is constant on connected-components of $\catI a$, or
equivalently that if $f_\blb: (a_\blb, m_\blb) \go (b_\blb, p_\blb)$ in
$\catI a$ then $\beta_a(a_\blb, m_\blb) = \beta_a(b_\blb, p_\blb)$.
Indeed, given such an $f_\blb$, there are coalgebra maps
\[
\begin{diagram}[width=6em]
(H^{(a_\sblb, m_\sblb)}, \theta^{(a_\sblb, m_\sblb)})	&
	&	&	&	\\
\uTo<{f_\sblb^*}	&
\rdTo(2,1)^{\kappa^{(a_\sblb, m_\sblb)}}		&
(\oba\catI, \iota)	&
\rTo^\beta	&
(J, \gamma)	\\
(H^{(b_\sblb, p_\sblb)}, \theta^{(b_\sblb, p_\sblb)})	&
\ruTo(2,1)_{\kappa^{(b_\sblb, p_\sblb)}}		&
	&	&	\\
\end{diagram}
\]
and $\beta \of \kappa^{(a_\sblb, m_\sblb)} \of f_\blb^* = \beta \of
\kappa^{(\beta_\sblb, p_\sblb)}$ by terminality of $(J, \gamma)$.  But
\begin{eqnarray*}
\kappa^{(a_\sblb, m_\sblb)} f_\blb^* (1_a)	&
=	&
\kappa^{(a_\sblb, m_\sblb)} (1_a \of f_0)
=
\kappa^{(a_\sblb, m_\sblb)} (1_a)
=
(a_\blb, m_\blb),
\\
\kappa^{(b_\sblb, p_\sblb)} (1_a)	&
=	&
(b_\blb, p_\blb),
\end{eqnarray*}
so $\beta_a (a_\blb, m_\blb) = \beta_a (b_\blb, p_\blb)$ as
required.  

The transformation $\ovln{\beta}: I \go J$ induces a functor $\elt{I} \go
\elt{J}$ over $\scat{A}$.  Moreover, there is a natural transformation
$\id_\Cat \go D \Pi_0$ (the unit of the adjunction $\Pi_0 \ladj D$), hence
a transformation $\catI \go D \Pi_0 \catI = D I$, hence a functor
$\elt{\catI} \go \elt{DI} \iso \elt{I}$ over $\scat{A}$.  So there is a
commutative triangle
\[
\begin{diagram}
\elt{\catI}	&		&\rTo^F	&		&\elt{J}	\\
		&\rdTo<\pr	&	&\ldTo>\pr	&		\\
		&		&\scat{A}&		&		\\
\end{diagram}
\]
in $\Cat$, where $\pr$ denotes a projection.  Now, condition \So\ says that
if $\scat{K}$ is either of the categories $\littlepullback$ or
$\littleequalizer$ and if $G: \scat{K} \go \elt{\catI}$ then $\pr \of G$
admits a cone.  But given such a $G$, nondegeneracy of $J$ implies that
$F\of G$ admits a cone, so $\pr \of G = \pr \of F \of G$ admits a cone, as
required.